\documentclass[12pt]{article}

\usepackage{amsmath,amssymb,amsthm}

\begin{document}

\newcommand{\End}{{\rm{End}\ts}}
\newcommand{\non}{\nonumber}
\newcommand{\wt}{\widetilde}
\newcommand{\wh}{\widehat}
\newcommand{\ot}{\otimes}
\newcommand{\la}{\lambda}
\newcommand{\al}{\alpha}
\newcommand{\be}{\beta}
\newcommand{\ga}{\gamma}
\newcommand{\de}{\delta^{}}
\newcommand{\om}{\omega^{}}
\newcommand{\hra}{\hookrightarrow}
\newcommand{\ve}{\varepsilon}
\newcommand{\ts}{\,}
\newcommand{\qin}{q^{-1}}
\newcommand{\tss}{\hspace{1pt}}
\newcommand{\U}{ {\rm U}}
\newcommand{\Y}{ {\rm Y}}
\newcommand{\C}{\mathbb{C}\tss}
\newcommand{\Z}{\mathbb{Z}}
\newcommand{\A}{\mathcal{A}}
\newcommand{\ZZ}{{\rm Z}}
\newcommand{\gl}{\mathfrak{gl}}
\newcommand{\oa}{\mathfrak{o}}
\newcommand{\spa}{\mathfrak{sp}}
\newcommand{\g}{\mathfrak{g}}
\newcommand{\ka}{\mathfrak{k}}
\newcommand{\p}{\mathfrak{p}}
\newcommand{\sll}{\mathfrak{sl}}
\newcommand{\agot}{\mathfrak{a}}
\newcommand{\qdet}{ {\rm qdet}\ts}
\newcommand{\sdet}{ {\rm sdet}\ts}
\newcommand{\sgn}{ {\rm sgn}}
\newcommand{\Sym}{\mathfrak S}

\renewcommand{\theequation}{\arabic{section}.\arabic{equation}}

\newtheorem{thm}{Theorem}[section]
\newtheorem{lem}[thm]{Lemma}
\newtheorem{prop}[thm]{Proposition}
\newtheorem{cor}[thm]{Corollary}

\theoremstyle{definition}
\newtheorem{defin}[thm]{Definition}
\newtheorem{example}[thm]{Example}

\theoremstyle{remark}
\newtheorem{remark}[thm]{Remark}

\newcommand{\bth}{\begin{thm}}
\renewcommand{\eth}{\end{thm}}
\newcommand{\bpr}{\begin{prop}}
\newcommand{\epr}{\end{prop}}
\newcommand{\ble}{\begin{lem}}
\newcommand{\ele}{\end{lem}}
\newcommand{\bco}{\begin{cor}}
\newcommand{\eco}{\end{cor}}
\newcommand{\bde}{\begin{defin}}
\newcommand{\ede}{\end{defin}}
\newcommand{\bex}{\begin{example}}
\newcommand{\eex}{\end{example}}
\newcommand{\bre}{\begin{remark}}
\newcommand{\ere}{\end{remark}}

\newcommand{\bal}{\begin{aligned}}
\newcommand{\eal}{\end{aligned}}
\newcommand{\beq}{\begin{equation}}
\newcommand{\eeq}{\end{equation}}
\newcommand{\ben}{\begin{equation*}}
\newcommand{\een}{\end{equation*}}

\newcommand{\bpf}{\begin{proof}}
\newcommand{\epf}{\end{proof}}

\def\beql#1{\begin{equation}\label{#1}}

\title{\Large\bf  Coideal subalgebras in quantum affine algebras}
\author{{A. I. Molev, \quad E. Ragoucy\quad and\quad P. Sorba}}

\date{} 

\maketitle

\begin{abstract}
We introduce two subalgebras
in the type $A$ quantum affine algebra which
are coideals with respect to the Hopf
algebra structure.
In the classical limit $q\to 1$
each subalgebra specializes
to the enveloping algebra $\U(\ka)$,
where $\ka$ is a fixed point subalgebra of the loop algebra
$\gl_N[\la,\la^{-1}]$ with respect to a natural involution
corresponding to the embedding
of the orthogonal or symplectic Lie algebra into $\gl_N$.
We also give an equivalent presentation of these coideal subalgebras
in terms of generators and defining relations which have the form
of reflection-type equations.
We provide evaluation homomorphisms from these algebras
to the twisted quantized enveloping algebras
introduced earlier by Gavrilik and Klimyk and by Noumi.
We also construct
an analog of the quantum determinant for each of the algebras and show
that its coefficients belong to the center of the algebra.
Their images under the evaluation homomorphism provide a family
of central elements of the corresponding twisted quantized enveloping algebra.
\vspace{5 mm}

\noindent
Preprint LAPTH-927/02

\end{abstract}


\vspace{13 mm}

\noindent
School of Mathematics and Statistics\newline
University of Sydney,
NSW 2006, Australia\newline
alexm@maths.usyd.edu.au

\vspace{7 mm}

\noindent
LAPTH, Chemin de Bellevue, BP 110\newline
F-74941 Annecy-le-Vieux cedex, France\newline
ragoucy@lapp.in2p3.fr

\vspace{7 mm}

\noindent
LAPTH, Chemin de Bellevue, BP 110\newline
F-74941 Annecy-le-Vieux cedex, France\newline
sorba@lapp.in2p3.fr

\newpage

\section{Introduction}\label{sec:int}
\setcounter{equation}{0}

For a simple Lie algebra $\g$ over $\C$ consider the corresponding
{\it quantized enveloping algebra\/} $\U_q(\g)$; see Drinfeld~\cite{d:ha},
Jimbo~\cite{j:qd}. If $\ka$ is a subalgebra of $\g$ then
$\U(\ka)$ is a Hopf subalgebra of $\U(\g)$. However, $\U_q(\ka)$, even when
it is defined, need not be isomorphic to a Hopf subalgebra of $\U_q(\g)$.
In the case where $(\g,\ka)$ is a classical symmetric pair the
{\it twisted quantized enveloping algebra\/} $\U^{\rm tw}_q(\ka)$
was introduced by Noumi~\cite{n:ms} (type $A$ pairs)
and by Noumi and Sugitani~\cite{ns:qs} (remaining classical types).
This is a subalgebra and a left coideal
of the Hopf algebra $\U_q(\g)$ which specializes to $\U(\ka)$
as $q\to 1$. The algebras $\U^{\rm tw}_q(\ka)$
play an important role in the theory of quantum symmetric spaces
developed in \cite{n:ms} and \cite{ns:qs}.
In particular, in the type $A$, which we are only concerned with
in this paper, there are two twisted quantized enveloping algebras
$\U^{\rm tw}_q(\oa_N)$ and $\U^{\rm tw}_q(\spa_{2n})$
corresponding to the symmetric pairs
\beql{symmp}
\bal
{\rm AI:}&  \qquad\qquad (\gl_N,\oa_N),\\
{\rm AII:}&  \qquad\qquad (\gl_{2n},\spa_{2n}),
\eal
\end{equation}
respectively.
It was also shown by Noumi~\cite{n:ms} that
the algebra $\U^{\rm tw}_q(\oa_N)$ coincides with the one introduced
earlier by Gavrilik and Klimyk~\cite{gk:qd}.
The algebra $\U^{\rm tw}_q(\oa_N)$ also appears as the symmetry
algebra for the $q$-oscillator representation of
the quantized enveloping algebra $\U_q(\spa_{2n})$;
see Noumi, Umeda and Wakayama~\cite{nuw:qd, nuw:dp}.

In Noumi's approach, the defining relations for the
quantized algebras can be written in the form of a reflection-type
equation. A constant solution of the reflection equation provides
an embedding of the twisted quantized enveloping
algebra into $\U_q(\gl_N)$.
The quantum homogeneous spaces corresponding to
the remaining series of the classical symmetric pairs of type $A$
\beql{symmp3}
{\rm AIII:}  \qquad\qquad (\gl_{N},\gl_{N-l}\oplus \gl_{\ts l})
\end{equation}
were studied
by Dijkhuizen, Noumi and Sugitani~\cite{dn:fq}, \cite{dns:ma}.
A one-parameter family of the constant solutions
of the appropriate reflection equation was
produced in \cite{dns:ma}, although no reflection-type
presentation of the subalgebras of type $\U^{\rm tw}_q(\ka)$ were
formally introduced.

A different description of the coideal subalgebras of $\U_q(\g)$ associated
with an arbitrary irreducible symmetric pair $(\g,\ka)$
was given by
Letzter~\cite{l:cs, l:qs}. The subalgebras
are presented by generators and explicit relations depending
on the Cartan matrix of $\g$; see \cite{l:qs}.
In particular, this work demonstrates the importance of the
coideal property: it makes the construction of the twisted quantized
algebras essentially unique.

Natural infinite-dimensional analogs of
the symmetric pairs are provided
by involutive subalgebras in the polynomial current Lie algebras
$\g[x]=\g\ot\C[x]$ or loop algebras $\g[\la,\la^{-1}]=\g\ot\C[\la,\la^{-1}]$.
Let $(\g,\ka)$ be a symmetric pair and
$\g=\ka\oplus \p$ be the decomposition determined by the involution
$\theta$ of $\g$. So, $\ka$ and $\p$ are the eigenspaces of $\theta$
with the eigenvalues $1$ and $-1$, respectively.
Then the twisted polynomial current Lie algebra $\g[x]^{\theta}$
can be defined by
\beql{twpolc}
\g[x]^{\theta}=\ka\oplus \p\tss x\oplus\ka\tss x^2 \oplus \p\tss x^3 \oplus\cdots,
\end{equation}
or, equivalently, it is the fixed point subalgebra of $\g[x]$ with respect
to the extension of $\theta$ given by
\beql{thetaext}
\theta: A\tss x^p\mapsto (-1)^p\ts\theta(A)\tss x^p,\qquad A\in \g.
\end{equation}
As it was demonstrated by Drinfeld~\cite{d:ha}, the enveloping algebra
$\U(\g[x])$ admits a canonical deformation in the class of Hopf algebras.
The corresponding ``quantum" algebra is called the {\it Yangian\/}
and denoted by $\Y(\g)$.
For the case where $\theta$ is an involution of $\gl_N$
corresponding to the pair of type AI or AII,
quantum analogs of the symmetric pairs $(\gl_N[x],\gl_N[x]^{\theta})$
are provided by the Olshanski {\it twisted Yangians\/} \cite{o:ty}.
These are
coideal subalgebras in the Yangian $\Y(\gl_N)$ and each of them is a deformation
of the enveloping algebra $\U(\gl_N[x]^{\theta})$.

In the AIII case, the quantum analogs of the pairs $(\gl_N[x],\gl_N[x]^{\theta})$
are provided by the {\it reflection algebras\/} $B(N,l)$ which are
coideal subalgebras  of $\Y(\gl_N)$
originally introduced by Sklyanin~\cite{s:bc}.
Recently, these algebras and their representations were studied
in connection with the {\it NLS model\/};
see Liguori, Mintchev and Zhao~\cite{lmz:be},
Mintchev, Ragoucy and Sorba~\cite{mrs:ss},
Molev and Ragoucy~\cite{mr:rr}.

For the symmetric pairs $(\g[x],\g[x]^{\theta})$ of general types
the corresponding coideal subalgebras in the Yangian $\Y(\g)$ were recently
introduced by Delius, MacKay and Short~\cite{dms:br} in relation with
the principal chiral models with boundaries.
These subalgebras are given
in terms of the $Q$-presentation
of the Yangian.
A different
$R$-matrix presentation of coideal subalgebras in the (super) Yangian $\Y(\g)$
is given by Arnaudon, Avan, Cramp\'e, Frappat and Ragoucy
\cite{aacfr:rp}.
Field theoretical applications of the coideal subalgebras
in the quantum affine algebras have been studied
in a recent paper by
Delius and MacKay~\cite{dm:qg}.

In the case of the loop algebra $\wh{\g}=\g[\la,\la^{-1}]$ there is
another natural way (cf. \eqref{thetaext}) to extend
the involution $\theta$ of $\g$,
\beql{thetaextl}
\theta:   A\tss \la^p\mapsto \theta(A)\tss \la^{-p},\qquad A\in \g
\end{equation}
and thus to define the fixed point subalgebra $\wh{\g}^{\ts\theta}$.
In this paper we introduce certain quantizations of the
symmetric pairs $(\wh{\gl}_N, \wh{\gl}_N^{\ts\theta})$
associated with the involution $\theta$ corresponding
to the pairs AI and AII.
We define the {\it twisted $q$-Yangians\/}
$\Y^{\rm tw}_q(\oa_N)$ and $\Y^{\rm tw}_q(\spa_{2n})$
as subalgebras of the quantum affine algebra
$\U_q(\wh{\gl}_N)$.
They are left coideals with respect to the coproduct on
$\U_q(\wh{\gl}_N)$ and
specialize to $\U(\wh{\gl}_N^{\ts\theta})$ as $q\to 1$.

At this point we consider it necessary to comment on the terminology.
Although, as we have mentioned above, the Lie algebra $\wh{\gl}_N^{\ts\theta}$
is not a ``twisted" quantum affine algebra in the usual meaning,
we believe the names we use for the coideal subalgebras
can be justified having in mind their analogy with both
the twisted Yangians and the $q$-{\it Yangian\/}; cf. \cite{nt:yg}.
The latter is a subalgebra
of the quantum affine algebra $\U_q(\wh{\gl}_N)$ which can be regarded
as a $q$-analog of the usual Yangian $\Y(\gl_N)$; see also Section~\ref{sec:csaff}
below.

Our first main result is a construction of the evaluation homomorphisms
\beql{evtw}
\Y^{\rm tw}_q(\oa_N)\to \U^{\rm tw}_q({\oa}_N),
\qquad
\Y^{\rm tw}_q(\spa_{2n})\to \U^{\rm tw}_q({\spa}_{2n})
\end{equation}
to the corresponding twisted
quantized enveloping algebras of \cite{gk:qd} and \cite{n:ms}.
Note that an evaluation homomorphism
$
\U_q(\tss\wh{\g}\tss)\to \U_q(\g)
$
from the quantum affine algebra to the corresponding
quantized enveloping algebra
only exists if $\g$ is of $A$ type, and the same holds for
the case of the Yangians; see Jimbo~\cite{j:qu}, Drinfeld~\cite{d:ha}.
In both the cases the evaluation homomorphisms play an important role
in the representation theory of the quantum algebras;
see Chari and Pressley~\cite{cp:gq}.
An evaluation homomorphism from the twisted Yangian to the
corresponding enveloping algebra $\U(\oa_N)$ or $\U(\spa_{2n})$
does exist (see \cite{o:ty}, \cite{mno:yc})
and has many applications in the classical representation
theory; see e.g. \cite{m:ya} for an overview.
Note also that the existence of the homomorphisms \eqref{evtw} is not
directly related with the corresponding fact for the $A$ type algebras
but is quite a nontrivial property of the reflection equations
satisfied by the generators of the twisted $q$-Yangians.

Next we construct an analog of the quantum determinant
for each twisted $q$-Yangian and show that its coefficients belong
to the center of this algebra. The application of the evaluation
homomorphism \eqref{evtw} yields a family of central elements
in $\U^{\rm tw}_q({\oa}_N)$ and $\U^{\rm tw}_q({\spa}_{2n})$.
In the orthogonal case we also produce a `short' determinant-like formula
for this analog which employs a certain map from the symmetric
group into itself. This same map was used in the short formulas
for the Sklyanin determinants for the twisted Yangians; see \cite{m:ya}.
Some other families of Casimir elements
were constructed by Noumi, Umeda and Wakayama \cite{nuw:dp}
and by Gavrilik and Iorgov~\cite{gi:ce}.
It would be interesting to understand the relationship between
the families, as well as to investigate possible applications
to the study of the quantum Howe dual
pairs; cf. \cite{nuw:qd, nuw:dp}.

Another intriguing problem is to construct coideal subalgebras
of $\U_q(\wh{\gl}_N)$ of type AIII, i.e., to find
$q$-analogs of the Sklyanin reflection algebras $B(N,l)$ mentioned above.

\medskip

We are grateful to Gustav Delius, Masatoshi Noumi and T\^oru Umeda
for valuable discussions.
The financial support of the Australian
Research Council and the
{\it Laboratoire d'Annecy-le-Vieux de Physique Th\'eorique\/}
is acknowledged.

\section{Coideal subalgebras of $\U_q(\gl_N)$}\label{sec:cs}
\setcounter{equation}{0}

We shall use an $R$-matrix presentation of the algebra $\U_q(\gl_N)$.
Our main references are Jimbo~\cite{j:qu} and
Reshetikhin, Takhtajan and Faddeev~\cite{rtf:ql}.
We fix a complex parameter $q$ which is nonzero and not a root of unity.
Consider the $R$-matrix
\beql{rmatrixc}
R=q\ts\sum_i E_{ii}\ot E_{ii}+\sum_{i\ne j} E_{ii}\ot E_{jj}+
(q-\qin)\sum_{i<j}E_{ij}\ot E_{ji}
\end{equation}
which is an element of $\End\C^N\ot \End\C^N$, where
the $E_{ij}$ denote the standard matrix units and the indices run over
the set $\{1,\dots,N\}$. The $R$-matrix satisfies the Yang--Baxter equation
\beql{YBEconst}
R_{12}\ts R_{13}\ts  R_{23} =  R_{23}\ts  R_{13}\ts  R_{12},
\end{equation}
where both sides take values in $\End\C^N\ot \End\C^N\ot \End\C^N$ and
the subindices indicate the copies of $\End\C^N$, e.g.,
$R_{12}=R\ot 1$ etc.

The {\it quantized enveloping algebra\/} $\U_q(\gl_N)$ is generated
by elements $t_{ij}$ and $\bar t_{ij}$ with $1\leq i,j\leq N$
subject to the relations\footnote{Our $T$ and $\overline T$ correspond
to the $L$-operators $L^-$ and $L^+$, respectively, in the notation of \cite{rtf:ql}.}
\beql{defrel}
\bal
t_{ij}&=\bar t_{ji}=0, \qquad 1 \leq i<j\leq N,\\
t_{ii}\ts \bar t_{ii}&=\bar t_{ii}\ts t_{ii}=1,\qquad 1\leq i\leq N,\\
R\ts T_1T_2&=T_2T_1R,\qquad R\ts \overline T_1\overline T_2=
\overline T_2\overline T_1R,\qquad
R\ts \overline T_1T_2=T_2\overline T_1R.
\eal
\end{equation}
Here $T$ and $\overline T$ are the matrices
\beql{matrt}
T=\sum_{i,j}t_{ij}\ot E_{ij},\qquad \overline T=\sum_{i,j}
\overline t_{ij}\ot E_{ij},
\end{equation}
which are regarded as elements of the algebra $\U_q(\gl_N)\ot \End\C^N$.
Both sides of each of the $R$-matrix relations in \eqref{defrel}
are elements of $\U_q(\gl_N)\ot \End\C^N\ot \End\C^N$ and the subindices
of $T$ and $\overline T$ indicate the copies of $\End\C^N$ where
$T$ or $\overline T$ acts; e.g. $T_1=T\ot 1$. In terms of the
generators the defining relations between the $t_{ij}$
can be written as
\beql{defrelg}
q^{\delta_{ij}}\ts t_{ia}\ts t_{jb}-
q^{\delta_{ab}}\ts t_{jb}\ts t_{ia}
=(q-\qin)\ts (\de_{b<a} -\de_{i<j})
\ts t_{ja}\ts t_{ib},
\end{equation}
where $\de_{i<j}$ equals $1$ if $i<j$ and $0$ otherwise.
The relations between the $\bar t_{ij}$
are obtained by replacing $t_{ij}$ by $\bar t_{ij}$ everywhere in
\eqref{defrelg}. Finally, the relations involving both
$t_{ij}$ and $\bar t_{ij}$ have the form
\beql{defrelg2}
q^{\delta_{ij}}\ts \bar t_{ia}\ts t_{jb}-
q^{\delta_{ab}}\ts t_{jb}\ts \bar t_{ia}
=(q-\qin)\ts (\de_{b<a}\ts t_{ja}\ts \bar t_{ib} -\de_{i<j}\ts
\ts \bar t_{ja}\ts t_{ib}).
\end{equation}

We shall also use another $R$-matrix $\wt R$ given by
\beql{rmatrixtld}
\wt R=\qin\ts\sum_i E_{ii}\ot E_{ii}+\sum_{i\ne j} E_{ii}\ot E_{jj}+
(\qin-q)\sum_{i>j}E_{ij}\ot E_{ji}.
\end{equation}
We have the relations
\beql{relRRt}
R-\wt R=(q-\qin)\ts P,\qquad \wt R=PR^{-1}P,
\end{equation}
where
\beql{p}
P=\sum_{i,j}E_{ij}\ot E_{ji}
\end{equation}
is the permutation operator.
The following relations are implied by \eqref{defrel}
\beql{relco}
\wt R\ts T_1T_2=T_2T_1\wt R,\qquad \wt R\ts \overline T_1\overline T_2=
\overline T_2\overline T_1\wt R,\qquad
\wt R\ts T_1\overline T_2=\overline T_2T_1\wt R.
\end{equation}

The coproduct $\Delta$ on $\U_q(\gl_N)$ is defined by
the relations
\beql{copr}
\Delta(t_{ij})=\sum_{k=1}^N t_{ik}\ot t_{kj},\qquad
\Delta(\bar t_{ij})=\sum_{k=1}^N \bar t_{ik}\ot \bar t_{kj}.
\end{equation}

It is well known that the algebra $\U_q(\gl_N)$
specializes to $\U(\gl_N)$  as $q\to 1$.
To make this more precise, regard $q$ as a formal variable
and $\U_q(\gl_N)$ as an algebra over $\C(q)$.
Then set $\A=\C[q,\qin]$ and consider
the $\A$-subalgebra $\U_{\A}$ of $\U_q(\gl_N)$ generated by
the elements
\beql{taugen1}
\tau_{ij}=\frac{t_{ij}}{q-\qin}\quad\text{for}\quad i> j,
\qquad
\bar\tau_{ij}=\frac{\bar t_{ij}}{q-\qin}\quad\text{for}\quad i< j,
\end{equation}
and
\beql{taugen2}
\tau_{ii}=\frac{t_{ii}-1}{q-1},\qquad
\bar\tau_{ii}=\frac{\bar t_{ii}-1}{q-1},
\end{equation}
for $i=1,\dots,N$. Then we have an isomorphism \beql{isom}
\U_{\A}\ot_{\A}\C\cong \U(\gl_N) \end{equation} with the action of $\A$ on
$\C$ defined via the evaluation $q=1$; see e.g.
\cite[Section~9.2]{cp:gq}. Note that $\tau_{ij}$ and
$\bar\tau_{ij}$ respectively specialize to the elements $E_{ij}$
and $-E_{ij}$ of $\U(\gl_N)$. More generally, given a subalgebra
$V$ of $\U_q(\gl_N)$, set $V_{\A}=V\cap \U_{\A}$. Following
Letzter~\cite[Section~1]{l:cs}, we shall say that $V$ specializes
to the subalgebra $V^{\circ}$ of $\U(\gl_N)$ (as $q$ goes to $1$)
if the image of $V_{\A}$ in $\U_{\A}\ot_{\A}\C$ is $V^{\circ}$.

\subsection{Orthogonal case}\label{subsec:oc}

Following Noumi~\cite{n:ms}, we introduce the {\it twisted
quantized enveloping algebra\/} $\U^{\rm tw}_q(\oa_N)$
as the subalgebra of $\U_q(\gl_N)$ generated by the matrix
elements of the matrix $S=T\ts \overline T^{\ts t}$. It can
be easily derived from \eqref{defrel} (see \cite{n:ms})
that the matrix $S$ satisfies the relations
\begin{align}\label{sijo}
s_{ij}&=0, \qquad 1 \leq i<j\leq N,\\
\label{sii1}
s_{ii}&=1,\qquad 1\leq i\leq N,\\
\label{rsrs}
R\ts S_1& R^{\ts t} S_2=S_2R^{\ts t} S_1R,
\end{align}
where $R^{\ts t}:=R^{\ts t_1}$ denotes the element obtained from $R$ by
the transposition in the first tensor factor:
\beql{rt}
R^{\ts t}=q\ts\sum_i E_{ii}\ot E_{ii}+\sum_{i\ne j} E_{ii}\ot E_{jj}+
(q-\qin)\sum_{i<j}E_{ji}\ot E_{ji}.
\end{equation}
Indeed, the only nontrivial part of this derivation is to verify
that
\beql{ssubt}
R\ts T_1\ts\overline T_1^t R^{\ts t}\ts T_2\ts\overline T_2^t=
T_2\ts\overline T_2^t R^{\ts t}\ts T_1\ts\overline T_1^t R.
\end{equation}
However, this is implied by the relation $R \ts R^{\ts t}=R^{\ts t} R$
and the following consequences of
\eqref{defrel}:
\beql{trt}
\overline T_1^t R^{\ts t}\ts T_2=T_2     R^{\ts t}\ts \overline T_1^t,\qquad
R\ts  \overline T_1^t \overline T_2^t=  \overline T_2^t  \overline T_1^t R.
\end{equation}

We now prove an auxiliary lemma which establishes
a weak form of the Poincar\'e--Birkhoff--Witt theorem
for abstract algebras defined by the relation \eqref{rsrs}.
It will be used in both the orthogonal and symplectic cases.

\ble\label{lem:wpbw}
Consider the associative algebra with $N^2$ generators
$s_{ij}$, $i,j=1,\dots,N$ and the defining relations
written in terms of the matrix $S=(s_{ij})$ by the relation \eqref{rsrs}.
Then the ordered monomials of the form
\beql{monord}
s_{11}^{k_{11}}\ts s_{12}^{k_{12}}\cdots s_{1N}^{k_{1N}}\  \cdots\
s_{N1}^{k_{N1}}\ts s_{N2}^{k_{N2}}\cdots s_{NN}^{k_{NN}}
\end{equation}
with nonnegative powers $k_{ij}$ linearly span the algebra.
\ele

\bpf  Rewriting \eqref{rsrs} in terms of the generators we get
\beql{drabs}
\bal
q^{\delta_{aj}+\delta_{ij}}\ts s_{ia}\ts s_{jb}-
q^{\delta_{ab}+\delta_{ib}}\ts s_{jb}\ts s_{ia}
{}&=(q-\qin)\ts q^{\delta_{ai}}\ts (\de_{b<a} -\de_{i<j})
\ts s_{ja}\ts s_{ib}\\
{}&+(q-\qin)\ts \big(q^{\delta_{ab}}\ts \de_{b<i}\ts s_{ji}\ts s_{ba}
- q^{\delta_{ij}}\ts \de_{a<j}\ts s_{ij}\ts s_{ab}\big)\\
{}&+ (q-\qin)^2\ts  (\de_{b<a<i} -\de_{a<i<j})\ts s_{ji}\ts s_{ab},
\eal
\end{equation}
where $\de_{i<j}$ or $\de_{i<j<k}$ equals $1$ if the subindex inequality
is satisfied and $0$ otherwise.
We shall be proving that any monomial
in the generators can be written as a linear combination of
monomials of the form \eqref{monord}. Given a monomial
\beql{anymon}
s_{i_1a_1}\cdots\ts s_{i_pa_p}
\end{equation}
we introduce its {\it length\/}
$p$ and {\it weight\/} $w$ by $w=i_1+\cdots+i_p$. Clearly,
for a monomial of length $p$ the weight can range between $p$ and $pN$.
We shall use induction on $w$.
By \eqref{drabs} we obtain the following
equality modulo products of weight less than $i+j$: for $i\geq j$
\beql{igj}
q^{\delta_{aj}+\delta_{ij}}\ts s_{ia}\ts s_{jb}\equiv
q^{\delta_{ab}+\delta_{ib}}\ts s_{jb}\ts s_{ia}
+\de_{b<a}\ts(q-\qin)\ts q^{\delta_{ai}}
\ts s_{ja}\ts s_{ib}.
\end{equation}
This allows us to represent \eqref{anymon} modulo
monomials of weight less than $w$
as a linear combination of monomials
$s_{j_1b_1}\cdots\ts s_{j_pb_p}$ of weight $w$ such that
$j_1\leq\cdots \leq j_p$. Consider a submonomial
$s_{ic_1}\cdots\ts s_{ic_r}$
containing generators with the same first index.
By \eqref{igj} we have for $a>b$
\beql{igi}
q^{\delta_{ai}+1}\ts s_{ia}\ts s_{ib}\equiv
q^{\delta_{ab}+\delta_{ib}}\ts s_{ib}\ts s_{ia}.
\end{equation}
Using this relation repeatedly we
bring the submonomial to the required
form.
\epf

We shall be proving now
that \eqref{sijo}--\eqref{rsrs} are precisely the defining relations of
the algebra $\U^{\rm tw}_q(\oa_N)$. In other words, the following theorem
takes place.

\bth\label{thm:inj}
The abstract algebra $\cal S$
generated by elements $s_{ij}$, $i,j=1,\dots,N$
with the defining relations \eqref{sijo}--\eqref{rsrs} is isomorphic to
$\U^{\rm tw}_q(\oa_N)$.
\eth

\bpf
It was mentioned in
\cite[Remark~7.9\tss(1)]{nuw:dp} without a detailed proof
that this fact can be established
with the use of the Diamond Lemma. We employ a different approach based
on a weak Poincar\'e--Birkhoff--Witt theorem for the algebra $\cal S$.

For this proof only, denote the matrix $T\ts\overline T^t$ by
$\wt S$ and its matrix elements by $\wt s_{ij}$.
As we noted above,
the map $S\mapsto \wt S$ defines an algebra homomorphism
$\varphi:{\cal S}\to\U^{\rm tw}_q(\oa_N)$. We only need to show that this
homomorphism is injective. By Lemma~\ref{lem:wpbw} the monomials
\beql{ordmon}
{s}_{21}^{\ts k_{21}}\ts    {s}_{31}^{\ts k_{31}}\ts
{s}_{32}^{\ts k_{32}}\ts  \cdots \ts
{s}_{N1}^{\ts k_{N1}}\ts {s}_{N2}^{\ts k_{N2}} \ts
\cdots\ts  {s}_{N,N-1}^{\ts k_{N,N-1}}
\end{equation}
span the algebra $\cal S$. We shall show that the images
of these monomials under $\varphi$ are linearly independent.
We prove, in fact, that given any linear ordering on the
set of generators $\wt s_{ij}$, $i>j$,
the ordered monomials in the
$\wt s_{ij}$ are linearly independent.
Regarding $q$ as a formal variable, we keep the notation
$\A$ for the algebra of Laurent polynomials $\C[q,\qin]$; see
Section~\ref{sec:cs}. Set $V=\U^{\rm tw}_q(\oa_N)$ and note that
the subalgebra $V_{\A}$ is generated by the elements
$\sigma_{ij}:=\wt s_{ij}/(q-\qin)$ with $i>j$.
It is enough
to verify that the ordered monomials in the generators $\sigma_{ij}$
are linearly independent over $\A$.
We have for $i>j$
\beql{sigmatau}
\sigma_{ij}=\tau_{ij}+\bar\tau_{ji}+(q-1)(\tau_{ii}\bar\tau_{ji}+
\tau_{ij}\bar\tau_{jj})+(q-\qin)\sum_{j<a<i}
\tau_{ia}\bar\tau_{ja}.
\end{equation}
Therefore, the image of $\sigma_{ij}$ in $V_{\A}\ot_{\A}\C$ is
$F_{ij}:=E_{ij}-E_{ji}$. The elements $F_{ij}$ with $i>j$
constitute a basis of a subalgebra of $\gl_N$ isomorphic to
the orthogonal Lie algebra $\oa_N$.
By the Poincar\'e--Birkhoff--Witt theorem for $\oa_N$
the ordered monomials in the $F_{ij}$ are linearly independent.
Now suppose, on the contrary, that there exists a nontrivial linear combination
of the ordered monomials in the $\sigma_{ij}$ equal to zero:
\beql{lincomb}
\sum_{(k)}c_{(k)}\prod_{i>j} \sigma_{ij}^{k_{ij}}=0,
\end{equation}
summed over the multi-indices $(k)=(k_{ij}\ts|\ts i>j)$, where $c_{(k)}\in\A$.
We may assume that
at least one coefficient $c_{(k)}$ in \eqref{lincomb} does not
vanish at $q=1$.
Taking the image of \eqref{lincomb} in $V_{\A}\ot_{\A}\C$ yields a nontrivial
linear combination of the ordered monomials in the $F_{ij}$ equal to zero.
This makes a contradiction proving the claim.
\epf

As a corollary of the above argument we obtain an analog of the
Poincar\'e--Birkhoff--Witt theorem for the algebra $\cal S$.
A different proof was given by N.~Iorgov; see~\cite{hkp:ce}.

\bco\label{cor:pbw}
The monomials \eqref{ordmon} constitute a basis of the algebra $\cal S$.
\qed
\eco

It follows from the proof of the theorem that the algebra $V_{\A}\ot_{\A}\C$
coincides with the enveloping algebra $\U(\oa_N)$.
Thus, $\U^{\rm tw}_q(\oa_N)$ specializes to
$\U(\oa_N)$ as $q\to 1$; see \cite{gk:qd} and \cite[Section~2.4]{n:ms}.
This result also follows from \cite[Section~6]{l:sp}.

Regarding $\U^{\rm tw}_q(\oa_N)$ as a subalgebra of $\U_q(\gl_N)$
we introduce another matrix $\overline S$ by
\beql{barS}
\overline S=\overline T\ts T^t.
\end{equation}
The matrix $\overline S=(\bar s_{ij})$ is upper triangular with 1's on the diagonal.
It is related to $S$ by the formula
\beql{SSbar}
\overline S=1-q +q\ts S^t,
\end{equation}
or, in terms of the matrix elements, $\bar s_{ij}=q\ts s_{ji}$
for $i<j$. Indeed, \eqref{defrel} implies
\beql{tbart}
\bar t_{ia}\ts t_{ja}=q\ts t_{ja}\ts \bar t_{ia}.
\end{equation}
Taking the sum over $a$ gives the result.
So, the elements $\bar s_{ij}$ belong to the subalgebra
$\U^{\rm tw}_q(\oa_N)$.

The next proposition is immediate from \eqref{copr}.
It shows that
the subalgebra $\U^{\rm tw}_q(\oa_N)$ is a left coideal of $\U_q(\gl_N)$.
This property was mentioned in \cite[Section~2.4]{n:ms} and
proofs can be found in \cite{gik:nd} and \cite[Lemma~6.3]{l:sp}.

\bpr\label{prop:coid}
The image of the generator $s_{ij}$ of $\U^{\rm tw}_q(\oa_N)$
under the coproduct
is given by
\beql{coprs}
\Delta(s_{ij})=\sum_{k,l=1}^N t_{ik}\ts\bar t_{jl}\ot s_{kl}.
\end{equation}
\qed
\epr

\bre
It was shown in \cite{n:ms} that
for any nondegenerate diagonal matrix $D$ the matrix
$S=T\tss D\tss \overline T^{\ts t}$ satisfies the
reflection equation \eqref{rsrs}. Therefore one can define
a family of subalgebras in $\U_q(\gl_N)$ parametrized
by the matrices $D$. However, all of them are isomorphic to each
other as abstract algebras which can be seen from the fact that
for any diagonal matrix $C$ the relation \eqref{rsrs}
is preserved by the transformation
$S\mapsto C\tss S\tss C$. Indeed, the entries of $S$ are then
transformed as $s_{ij}\mapsto s_{ij}\ts c_i\ts c_j$ and the claim
is immediate from \eqref{drabs}.
The corresponding remark also applies
to the symplectic case
considered below.
\ere

\subsection{Symplectic case}\label{subsec:sc}

Following again \cite{n:ms} introduce the $2n\times 2n$ matrix $G$ by
\beql{g}
G=q\ts \sum_{k=1}^n E_{2k-1,2k}-\sum_{k=1}^n E_{2k,2k-1}.
\end{equation}
It is convenient to introduce the involution of the set
of indices $\{1,\dots,2n\}$ which we denote by a prime
and which acts by the rule: for $k=1,\dots,n$
\beql{invols}
(2k-1)'=2k\qquad\text{and}\qquad (2k)'=2k-1.
\end{equation}
To introduce
the
{\it twisted
quantized enveloping algebra\/} $\U^{\rm tw}_q(\spa_{2n})$
consider the matrix $S=T\ts G\ts \overline T^{\ts t}$
with entries in $\U_q(\gl_{2n})$. For each odd $i=1,3,\dots,2n-1$
the element $s_{ii'}=q\ts t_{ii}\ts\bar t_{i'i'}$ is invertible
in $\U_q(\gl_{2n})$. We define $\U^{\rm tw}_q(\spa_{2n})$
as the subalgebra of $\U_q(\gl_{2n})$
generated by the matrix
elements of the matrix $S$ and by
the elements\footnote{These additional generators,
not considered in \cite{n:ms}, will bring some simplification
in the arguments below.} $s_{ii'}^{-1}$,
$i=1,3,\dots,2n-1$.
Clearly, the matrix $S$ has a block-triangular form with
$n$ diagonal $2\times 2$-blocks so that
\beql{szer}
s_{ij}=0\qquad\text{for}\quad i<j \quad\text{with}\quad j\ne i'.
\end{equation}
As in the orthogonal case, the matrix $S$ satisfies the relation
\beql{refsymp}
R\ts S_1 R^{\ts t} S_2=S_2R^{\ts t} S_1R,
\end{equation}
with the same definition of $R$ and $R^{\ts t}$ (taking $N=2n$),
cf. \eqref{sijo}--\eqref{rsrs}.
This follows from the fact that
$G$ is a solution of the reflection equation \cite{n:ms}
\beql{refg}
R\ts G_1 R^{\ts t} G_2=G_2R^{\ts t} G_1R.
\end{equation}

\ble\label{lem:ker}
For any odd $i=1,3,\dots,2n-1$ we have the identity
\beql{qds}
s_{i'i'}\ts s_{ii}-q^2\ts s_{i'i}\ts s_{ii'}=q^3.
\end{equation}
\ele

\bpf The calculation is the same for each $i$ so we take $i=1$.
We have
\beql{stt}
\bal
s_{11}&=q\ts t_{11}\ts\bar t_{12},\qquad s_{12}=q\ts t_{11}\ts\bar t_{22}\\
s_{22}&=q\ts t_{21}\ts\bar t_{22},\qquad s_{21}=q\ts t_{21}\ts\bar t_{12}
- t_{22}\ts\bar t_{11}.
\eal
\end{equation}
Now the relation is implied by the defining relations \eqref{defrel}
in $\U_q(\gl_N)$.
\epf

Next we prove an analog of Theorem~\ref{thm:inj} for the algebra
$\U^{\rm tw}_q(\spa_{2n})$. Introduce the abstract algebra $\cal S$
with generators $s_{ij}$, $i,j=1,2,\dots,2n$
and $s_{ii'}^{-1}$, $i=1,3,\dots,2n-1$ with the defining
relations given by \eqref{szer}, \eqref{refg}, \eqref{qds} and
\beql{ssinv}
s_{ii'}\ts s_{ii'}^{-1}=s_{ii'}^{-1}\ts s_{ii'}=1,\qquad
i=1,3,\dots,2n-1.
\end{equation}
We show first that the relations \eqref{qds} allow us to eliminate
the generators $s_{i'i}$ in the spanning set of monomials.

\ble\label{lem:bass}
The algebra $\cal S$ is spanned by the
ordered monomials of the form
\beql{ordsymp}
\prod_{i=1,3,\dots,2n-1}^{\rightarrow}
s_{i1}^{\ts k_{i1}}\ts s_{i2}^{\ts k_{i2}} \ts
\cdots\ts  s_{ii'}^{\ts k_{ii'}}\ts
s_{i'i'}^{\ts k_{i'i'}}\ts s_{i'1}^{\ts k_{i'1}} \ts
\cdots\ts  s_{i',i'-2}^{\ts k_{i',i'-2}},
\end{equation}
where the $k_{ii'}$ with $i=1,3,\dots,2n-1$ are arbitrary integers while
the remaining powers $k_{ij}$ are nonnegative integers.
\ele

\bpf Denote by $\cal S'$ the algebra with generators $s_{ij}$
and $s_{ii'}^{-1}$ and the defining
relations given by \eqref{szer}, \eqref{refg} and
\eqref{ssinv}. The defining relations imply that
for any odd $i$ we have
\beql{sii+1}
q^{\delta_{i',k}-\delta_{ik}} s_{ii'}s_{kl} =
q^{\delta_{il}-\delta_{i',l}} s_{kl}  s_{ii'}.
\end{equation}
Therefore, by Lemma~\ref{lem:wpbw}
the algebra $\cal S'$ is spanned by the monomials
\beql{ordsympf}
\prod_{i=1,3,\dots,2n-1}^{\rightarrow}
s_{i1}^{\ts k_{i1}}\ts s_{i2}^{\ts k_{i2}} \ts
\cdots\ts  s_{ii'}^{\ts k_{ii'}}\ts
 s_{i'1}^{\ts k_{i'1}} \ts
\cdots\ts  s_{i',i'-1}^{\ts k_{i',i'-1}}\ts s_{i'i'}^{\ts k_{i'i'}},
\end{equation}
where the integers $k_{ii'}$ are allowed to be negative.
As we noted in the proof of
Lemma~\ref{lem:wpbw}, the generators $s_{ia}$ and $s_{ib}$
can be permuted modulo terms of lower weight; see \eqref{igi}.
Therefore, rearranging the generators appropriately, we conclude that
$\cal S'$ is also spanned by the monomials
\beql{ordsympf2}
\prod_{i=1,3,\dots,2n-1}^{\rightarrow}
s_{i1}^{\ts k_{i1}}\ts s_{i2}^{\ts k_{i2}} \ts
\cdots\ts  s_{ii'}^{\ts k_{ii'}}\ts
s_{i',i'-1}^{\ts k_{i',i'-1}}\ts s_{i'i'}^{\ts k_{i'i'}}
\ts s_{i'1}^{\ts k_{i'1}} \ts
\cdots\ts  s_{i',i'-2}^{\ts k_{i',i'-2}}.
\end{equation}
It is a straightforward calculation to derive from the defining
relations that the elements $s_{i'i'}\ts s_{ii}-q^2\ts s_{i'i}\ts s_{ii'}$
with $i=1,3,\dots,2n-1$ are central in the algebra $\cal S'$.
Let $I$ be the ideal of $\cal S'$ generated by the central
elements $s_{i'i'}\ts s_{ii}-q^2\ts s_{i'i}\ts s_{ii'}-q^3$
for $i=1,3,\dots,2n-1$. We need to show that
the quotient ${\cal S}'/I$ is spanned by
the monomials \eqref{ordsymp}. However,
the defining relations between the generators
$s_{ii}$, $s_{ii'}$, $s_{i'i}$ and $s_{i'i'}$
do not involve any other generators. In order to complete the proof
it is therefore sufficient to consider the particular case $n=1$.
We shall show that
modulo the ideal $I$ generated by the central element
$s_{22}\ts s_{11}-q^2\ts s_{21}\ts s_{12}-q^3$
every element of the algebra $\cal S'$ can be written as a linear
combination of monomials of the form
$s_{11}^{\ts k_{11}}\ts s_{12}^{\ts k_{12}}\ts
{s}_{22}^{\ts k_{22}}$.
Let $s_{11}^{\ts k_{11}}\ts s_{12}^{\ts k_{12}} \ts
\ts s_{21}^{\ts k_{21}}\ts
{s}_{22}^{\ts k_{22}}$ be an arbitrary monomial of the form
\eqref{ordsympf2}. We assume $k_{21}\geq 1$ and use induction on $k_{21}$.
Modulo the ideal $I$ we have
\beql{s21eq}
s_{21}\equiv q^{-2}s_{22}\ts s_{11}\ts s_{12}^{-1}+q\ts s_{12}^{-1}.
\end{equation}
Therefore, omitting the ordered monomials with smaller powers $k_{21}$
we have the relation
modulo the ideal $I$:
\beql{eqmi}
s_{11}^{\ts k_{11}}\ts s_{12}^{\ts k_{12}}  \ts
\ts s_{21}^{\ts k_{21}}\ts
{s}_{22}^{\ts k_{22}}\equiv q^{-2}s_{11}^{\ts k_{11}}\ts s_{12}^{\ts k_{12}-1}
\ts s_{22}\ts s_{11}\ts s_{21}^{\ts k_{21}-1}\ts {s}_{22}^{\ts k_{22}}.
\end{equation}
Note that by \eqref{drabs} we have the relation
\beql{s22s11}
s_{22}\ts s_{11}=s_{11}\ts s_{22}+(q-\qin)(s_{12}^2+q\ts s_{12}\ts s_{21}).
\end{equation}
This allows us to bring the right hand side of \eqref{eqmi} to the form
\beql{ffmon}
(1-q^{-2})\ts s_{11}^{\ts k_{11}}\ts s_{12}^{\ts k_{12}}    \ts
\ts s_{21}^{\ts k_{21}}\ts
{s}_{22}^{\ts k_{22}},
\end{equation}
which completes the proof.
\epf

The following is an analog of Theorem~\ref{thm:inj} for the
symplectic case.

\bth\label{thm:isom} The algebra $\cal S$ is isomorphic to
$\U^{\rm tw}_q(\spa_{2n})$.
\eth

\bpf We use the same argument as for the proof of Theorem~\ref{thm:inj}.
For this proof only, denote the matrix $T\ts G\ts \overline T^{\ts t}$ by
$\wt S$ and its matrix elements by $\wt s_{ij}$.
The map $S\mapsto \wt S$ defines an algebra homomorphism
$\varphi:{\cal S}\to\U^{\rm tw}_q(\spa_{2n})$. We show that the images
of the monomials \eqref{ordsymp} under $\varphi$ are linearly independent.
Note that by
\eqref{sii+1} the product of
a monomial of the form \eqref{ordsymp} and $s_{ii'}^k$ is,
up to a nonzero factor,
equal to the same monomial with the index $k_{ii'}$
replaced with  $k_{ii'}+k$. Therefore we may assume without loss of
generality that
all powers in \eqref{ordsymp} are nonnegative integers.

Set $V=\U^{\rm tw}_q(\spa_{2n})$.
The corresponding generators $\sigma_{ij}\in V_{\A}$
can now be given by
\beql{sigsym}
\sigma_{ij}=\frac{\wt s_{ij}-g_{ij}}{q-\qin},
\end{equation}
where the $g_{ij}$ are the matrix elements of $G$.
We have the relation modulo $(q-1)$,
\beql{sigmataus}
\sigma_{ij}\equiv \wt g_{j'j}\ts\tau_{ij'}+\wt g_{ii'}\ts\bar\tau_{ji'},
\end{equation}
where $\wt g_{ij}$ is the value of $g_{ij}$ at $q=1$ so that
$\wt g_{ii'}=(-1)^{i-1}$.
Therefore, the image
of the element $\sigma_{ij}$ in $V_{\A}\ot_{\A}\C$ is
$F_{ij}:=\wt g_{j'j}\ts E_{ij'}-\wt g_{ii'}\ts E_{ji'}$.
The elements $F_{ij}$ span a subalgebra of $\gl_{2n}$ isomorphic to
the symplectic Lie algebra $\spa_{2n}$ and this allows us to complete
the proof exactly as in the orthogonal case.
\epf

The following is an analog of the
Poincar\'e--Birkhoff--Witt theorem for the algebra $\cal S$
which is immediate from Theorem~\ref{thm:isom}.

\bco\label{cor:pbws}
The monomials \eqref{ordsymp} constitute a basis of the algebra $\cal S$.
\qed
\eco

It follows from the proof of the theorem that the algebra $V_{\A}\ot_{\A}\C$
coincides with the enveloping algebra $\U(\spa_{2n})$.
Thus, $\U^{\rm tw}_q(\spa_{2n})$ specializes to
$\U(\spa_{2n})$ as $q\to 1$. The result also
follows from \cite[Section~6]{l:sp}.

Regarding $\U^{\rm tw}_q(\spa_{2n})$ as a subalgebra of $\U_q(\gl_{2n})$
we introduce another matrix $\overline S$ by
\beql{barSs}
\overline S=\overline T\ts G\ts T^t.
\end{equation}
Using the same argument as in the orthogonal case (see Section~\ref{subsec:oc})
we derive the following relations between the matrix
elements of $S$ and $\overline S$: for any $i=1,3,\dots 2n-1$
\beql{SSbars}
\bal
\bar s_{ii}&=-q^{-2}\ts s_{ii},\qquad \bar s_{i'i'}=-q^{-2}\ts s_{i'i'},\\
\bar s_{i'i}&=-q^{-1}\ts s_{ii'},\qquad
\bar s_{ii'}=-q^{-1}\ts s_{i'i}+(1-q^{-2})\ts s_{ii'},
\eal
\end{equation}
while for the remaining generators we have
\beql{remainij}
\bar s_{ij}=-q^{-1}\ts s_{ji}, \qquad i<j,\quad j\ne i'.
\end{equation}
Thus, the elements $\bar s_{ij}$ belong to the subalgebra
$\U^{\rm tw}_q(\spa_{2n})$.

The next proposition is immediate from \eqref{copr}; cf. Proposition~\ref{prop:coid}.
It shows that
the subalgebra $\U^{\rm tw}_q(\spa_{2n})$ is a left coideal of $\U_q(\gl_{2n})$;
see \cite[Lemma~6.3]{l:sp}.

\bpr\label{prop:coids}
The images of the generators of $\U^{\rm tw}_q(\spa_{2n})$
under the coproduct
are given by
\beql{coprss}
\Delta(s_{ij})=\sum_{k,l=1}^{2n} t_{ik}\ts\bar t_{jl}\ot s_{kl}
\end{equation}
and
\beql{coprsinv}
\Delta(s_{ii'}^{-1})=t_{i'i'}\ts\bar t_{ii}\ot s_{ii'}^{-1}.
\end{equation}
\qed
\epr

\section{Coideal subalgebras of $\U_q(\wh\gl_N)$}\label{sec:csaff}
\setcounter{equation}{0}

Consider the Lie algebra of Laurent polynomials $\gl_N[\la,\la^{-1}]$
in an indeterminate $\la$. We denote it by
$\wh\gl_N$ for brevity. The {\it quantum affine algebra\/}
$\U_q(\wh\gl_N)$ is a deformation of the universal enveloping algebra
$\U(\wh\gl_N)$. We use its $R$-matrix presentation following
Reshetikhin, Takhtajan and Faddeev~\cite{rtf:ql}; cf. Section~\ref{sec:cs}.
Note also a recent work of Frenkel and Mukhin~\cite{fm:ha},
where different realizations of $\U_q(\wh\gl_N)$ are collected.
By definition, the algebra $\U_q(\wh\gl_N)$ has countably many
generators $t_{ij}^{(r)}$ and $\bar t_{ij}^{\ts(r)}$ where
$1\leq i,j\leq N$ and $r$ runs over nonnegative integers.
They are combined into the matrices
\beql{taff}
T(u)=\sum_{i,j=1}^N t_{ij}(u)\ot E_{ij},\qquad
\overline T(u)=\sum_{i,j=1}^N \bar t_{ij}(u)\ot E_{ij},
\end{equation}
where $t_{ij}(u)$ and $\bar t_{ij}(u)$ are formal series
in $u^{-1}$ and $u$, respectively:
\beql{expa}
t_{ij}(u)=\sum_{r=0}^{\infty}t_{ij}^{(r)}\ts u^{-r},\qquad
\bar t_{ij}(u)=\sum_{r=0}^{\infty}\bar t_{ij}^{\ts(r)}\ts u^{r}.
\end{equation}
The defining relations are
\beql{defrelaff}
\bal
t_{ij}^{(0)}&=\bar t_{ji}^{\ts(0)}=0, \qquad 1 \leq i<j\leq N,\\
t_{ii}^{(0)}\ts \bar t_{ii}^{\ts(0)}&=\bar t_{ii}^{\ts(0)}
\ts t_{ii}^{(0)}=1,\qquad 1\leq i\leq N,\\
R(u,v)\ts T_1(u)T_2(v)&=T_2(v)T_1(u)R(u,v),\\
R(u,v)\ts \overline T_1(u)\overline T_2(v)&=
\overline T_2(v)\overline T_1(u)R(u,v),\\
R(u,v)\ts \overline T_1(u)T_2(v)&=T_2(v)\overline T_1(u)R(u,v),
\eal
\end{equation}
where we have used the notation of \eqref{defrel}
and $R(u,v)$ is
the trigonometric $R$-matrix given by
\beql{trRm}
\bal
R(u,v)={}&(u-v)\sum_{i\ne j}E_{ii}\ot E_{jj}+(\qin u-q\tss v)
\sum_{i}E_{ii}\ot E_{ii} \\
{}+ {}&(\qin-q)\tss u\tss\sum_{i> j}E_{ij}\ot
E_{ji}+ (\qin-q)\tss v\tss\sum_{i< j}E_{ij}\ot E_{ji}.
\eal
\end{equation}
It satisfies the Yang--Baxter equation
\beql{YBE}
R_{12}(u,v)  R_{13}(u,w)R_{23}(v,w) =  R_{23}(v,w) R_{13}(u,w) R_{12}(u,v),
\end{equation}
where both sides take values in $\End\C^N\ot \End\C^N\ot \End\C^N$ and
the subindices indicate the copies of $\End\C^N$, e.g.,
$R_{12}(u,v)=R(u,v)\ot 1$ etc. Note that $R(u,v)$ is related with
the constant $R$-matrices \eqref{rmatrixc} and \eqref{rmatrixtld} by
the formula
\beql{rwtR}
R(u,v)=u\ts \wt R-v\ts R.
\end{equation}

There is a Hopf algebra structure on $\U_q(\wh\gl_N)$ with the coproduct
defined by
\beql{copraff}
\Delta\big(t_{ij}(u)\big)=\sum_{k=1}^N t_{ik}(u)\ot t_{kj}(u),
\qquad
\Delta\big(\bar t_{ij}(u)\big)=\sum_{k=1}^N\bar t_{ik}(u)\ot\bar t_{kj}(u).
\end{equation}

The algebra $\U_q(\wh\gl_N)$ specializes to $\U(\wh\gl_N)$ as $q\to 1$.
More precisely,
as with the case of $\U_q(\gl_N)$ (see Section~\ref{sec:cs}),
regard $q$ as a formal variable
and introduce the $\A$-subalgebra $\U_{\A}$
of $\U_q(\wh\gl_N)$ generated by
the elements $\tau_{ij}^{(r)}$ and $\bar\tau_{ij}^{(r)}$
defined by
\beql{taugenaff1}
\tau_{ij}^{(r)}=\frac{t_{ij}^{(r)}}{q-\qin},\qquad
\bar\tau_{ij}^{(r)}=\frac{\bar t_{ij}^{\ts(r)}}{q-\qin}
\end{equation}
for $r\geq 0$ and all $i,j$, except for the case
$r=0$ and $i=j$ where we set
\beql{taugenaff2}
\tau_{ii}^{(0)}=\frac{t_{ii}^{\ts(0)}-1}{q-1},\qquad
\bar\tau_{ii}^{(0)}=\frac{\bar t_{ii}^{\ts(0)}-1}{q-1}.
\end{equation}
Then we have an isomorphism
\beql{limim}
\U_{\A}\ot_{\A}\C\cong\U(\wh\gl_N);
\end{equation}
see \cite[Section~12.2]{cp:gq} and \cite[Section~2]{fm:ha}.
The images of the generators of $\U_{\A}$ in \eqref{limim} are
given by
\beql{classlimaff}
\tau_{ij}^{(r)}\to E_{ij}\la^r,\qquad  \bar\tau_{ij}^{(r)}\to -E_{ij}\la^{-r}
\end{equation}
for all $r\geq 0$ with the exception $\tau_{ij}^{(0)}=\bar\tau_{ji}^{(0)}=0$
if $i<j$.
Given a subalgebra $V$ of $\U_q(\wh\gl_N)$ we set
$V_{\A}=V\cap \U_{\A}$. We shall say that $V$ specializes
to a subalgebra $V^{\circ}$ of $\U(\wh\gl_N)$ if
$V_{\A}\ot_{\A}\C\cong V^{\circ}$.

The quantized enveloping algebra
$\U_q(\gl_N)$ is a natural (Hopf) subalgebra of $\U_q(\wh\gl_N)$
defined by the embedding
\beql{emb}
t_{ij}\mapsto t_{ij}^{(0)},\qquad \bar t_{ij}\mapsto\bar t_{ij}^{\ts(0)}.
\end{equation}
Moreover, there is an algebra homomorphism $\U_q(\wh\gl_N)\to \U_q(\gl_N)$
called the {\it evaluation homomorphism\/} defined by
\beql{eval}
T(u)\mapsto T-\overline T\ts u^{-1},\qquad
\overline T(u)\mapsto \overline T-T\ts u.
\end{equation}
The $A$ type quantum affine algebras are exceptional
in the sense that only in this case such an evaluation homomorphism
does exist; see Chari--Pressley~\cite[Chapter~12]{cp:gq}.

The subalgebra of $\U_q(\wh\gl_N)$
generated by the elements $t_{ij}^{(r)}$
was studied e.g. in \cite{c:ni}, \cite{nt:yg}, \cite{rtf:ql}.
We call it the $q$-{\it Yangian\/}.
In what follows we construct quantum affine algebras associated with
the orthogonal and symplectic Lie algebras for which analogs of the
evaluation homomorphism \eqref{eval} do exist; cf. the $B$ and $C$
type twisted Yangians \cite{o:ty}, \cite{mno:yc}.
These algebras can be viewed as twisted analogs of
the $q$-Yangian as well as $q$-analogs
of the twisted Yangians. Note, however, that contrary to the case of the
twisted Yangians, our algebras are {\it not\/} subalgebras
of the $q$-Yangian; they are generated by certain combinations of
both types of elements $t_{ij}^{(r)}$ and $\bar t_{ij}^{\ts (r)}$.

\subsection{Orthogonal twisted $q$-Yangians}\label{subsec:oty}

\bde\label{def:orttwy} The {\em twisted $q$-Yangian\/} $\Y^{\rm tw}_q(\oa_{N})$
is the subalgebra of $\U_q(\wh\gl_N)$ generated by the coefficients $s_{ij}^{(r)}$ of
the matrix elements of the matrix $S(u)=T(u)\ts \overline T(u^{-1})^t$.
More precisely, we have
\beql{sumatr}
s_{ij}(u)=\sum_{a=1}^N t_{ia}(u)\ts \bar t_{ja}(u^{-1})
\end{equation}
so that
\beql{scoeff}
s_{ij}(u)=\sum_{r=0}^{\infty} s_{ij}^{(r)}\ts u^{-r}.
\end{equation}
The subalgebra $\Y^{\rm tw}_q(\oa_{N})$ is generated by
the elements $s_{ij}^{(r)}$ with $1\leq i,j\leq N$ and $r$ running over
the set of nonnegative integers.
\ede

Next we give a presentation of the algebra $\Y^{\rm tw}_q(\oa_{N})$
in terms of generators and defining relations by analogy
with the finite-dimensional case; see Section~\ref{subsec:oc}.
Consider the element $R^{\ts t}(u,v):=R^{\ts t_1}(u,v)$ obtained from $R(u,v)$
by the transposition in the first factor:
\beql{rtaff}
\bal
R^{\ts t}(u,v)={}&(u-v)\sum_{i\ne j}E_{ii}\ot E_{jj}+(\qin u-q\tss v)
\sum_{i}E_{ii}\ot E_{ii} \\
{}+ {}&(\qin-q)\tss u\tss\sum_{i> j}E_{ji}\ot
E_{ji}+ (\qin-q)\tss v\tss\sum_{i< j}E_{ji}\ot E_{ji}.
\eal
\end{equation}
The following relations are implied by \eqref{defrelaff}:
\begin{align} \label{sijoaff}
s_{ij}^{(0)}&=0, \qquad 1 \leq i<j\leq N,\\
\label{siio}
s_{ii}^{(0)}&=1, \qquad 1\leq i\leq N,\\
\label{rsrsaff}
R(u,v)\ts S_1(u)\ts R^{\ts t}(u^{-1},v)\ts S_2(v)&=S_2(v)
\ts R^{\ts t}(u^{-1},v)\ts S_1(u)\ts R(u,v).
\end{align}
We shall be proving that these are precisely the defining
relations for the algebra $\Y^{\rm tw}_q(\oa_{N})$.

\ble\label{lem:wpbwaff}
Consider the (abstract) associative algebra with generators
$s_{ij}^{(r)}$ where $i,j=1,\dots,N$ and $r=0,1,\dots$.
The defining relations are
written in terms of the matrix $S(u)=(s_{ij}(u))$ by the relation
\eqref{rsrsaff} with $s_{ij}(u)$ defined by \eqref{scoeff}.
Introduce the ordering on the generators in such a way
that $s_{ij}^{(r)}\preceq s_{kl}^{(p)}$ if and only if
$(r,i,j)\preceq (p,k,l)$
in the lexicographical order. Then the ordered monomials
in the generators span the algebra.
\ele

\bpf Write the defining relations in terms
of the generating series $s_{ij}(u)$:
\beql{draffsu}
\bal
(q^{-\delta_{ij}}u-q^{\delta_{ij}}v)
\ts\al_{ijab}(u,v)
+ (\qin-q)(u\de_{j<i}+v\de_{i<j})&
\ts\al_{jiab}(u,v)\\
=
(q^{-\delta_{ab}}u-q^{\delta_{ab}}v)
\ts\al_{jiba}(v,u)
+ (\qin-q)(u\de_{a<b}+v\de_{b<a})&
\ts\al_{jiab}(v,u),
\eal
\end{equation}
where we have used the notation
\beql{alphanot}
\al_{ijab}(u,v)=(q^{-\delta_{aj}}-q^{\delta_{aj}}uv)\ts s_{ia}(u)\ts s_{jb}(v)+
(\qin-q)(\de_{j<a}+uv\de_{a<j})\ts s_{ij}(u)\ts s_{ab}(v).
\end{equation}
We shall show that any monomial
\beql{mongene}
{s}_{i_1a_1}^{\ts (k_{1})}\ts {s}_{i_2a_2}^{\ts (k_{2})}
\cdots {s}_{i_pa_p}^{\ts (k_{p})}
\end{equation}
can be written as a linear combination of the ordered monomials.
Define the degree of the monomial \eqref{mongene} as the sum $k_1+\cdots+k_p$ and
argue by induction on the degree.
The induction base is Lemma~\ref{lem:wpbw} which takes care of the monomials
of degree zero. In other words, we can introduce the filtration on the algebra
by setting $\deg s_{ia}^{(k)}=k$ and it will be sufficient to prove the lemma
for the corresponding graded algebra. We keep the same notation
for its generators while the defining relations are given by \eqref{draffsu}
where instead of \eqref{alphanot} we should take
\beql{alphanotgr}
\al_{ijab}(u,v)=-q^{\delta_{aj}}\ts s_{ia}(u)\ts s_{jb}(v)+
(\qin-q)\ts\de_{a<j}\ts s_{ij}(u)\ts s_{ab}(v).
\end{equation}
Furthermore, we apply the same argument to this new algebra
considering the filtration defined by setting the degree of
$s_{ia}^{(k)}$ to be equal to $i$. The corresponding graded algebra
is generated by elements $s_{ia}^{(k)}$ with the defining relations
\eqref{draffsu} where the expression $\al_{ijab}(u,v)$ further simplifies
to
\beql{alphanotgrwt}
\al_{ijab}(u,v)=-q^{\delta_{aj}}\ts s_{ia}(u)\ts s_{jb}(v).
\end{equation}
Working with this algebra we have for $i>j$
\beql{reli>j}
\bal
q^{\delta_{aj}}\ts s_{ia}(u)\ts s_{jb}(v) &{}=\frac{q^{-\delta_{ab}}u-q^{\delta_{ab}}v}
{u-v}\ts  q^{\delta_{bi}}\ts s_{jb}(v)\ts  s_{ia}(u)\\
&{}+\frac{q-\qin}{u-v}\ts q^{\delta_{ai}}\Big(u\ts s_{ja}(u)\ts s_{ib}(v)
-(u\de_{a<b}+v\de_{b<a})\ts s_{ja}(v)\ts s_{ib}(u) \Big).
\eal
\end{equation}
This implies that the monomials \eqref{mongene} with the condition
$i_1\leq\cdots\leq i_p$ span the algebra.  Now we take $i=j$ and $a<b$ in
\eqref{draffsu} with the assumption \eqref{alphanotgrwt} to get
\beql{relijab}
\bal
q^{\delta_{bi}}\ts s_{ib}(v)\ts s_{ia}(u) &{}=
\frac{q^{-1}u-q\ts v}
{q^{-\delta_{ab}}u-q^{\delta_{ab}}v}\ts
q^{\delta_{ai}}\ts s_{ia}(u)\ts  s_{ib}(v)\\
&{}+\frac{q-\qin}{q^{-\delta_{ab}}u-q^{\delta_{ab}}v}\ts
q^{\delta_{ai}}\ts u\ts s_{ia}(v)\ts s_{ib}(u).
\eal
\end{equation}
Therefore, the algebra is spanned by the monomials \eqref{mongene} such that
$(i_1,a_1)\preceq\cdots\preceq (i_p,a_p)$ in the lexicographical order.
Finally, using again \eqref{draffsu} we note that
\beql{ijij}
s_{ia}(u)\ts s_{ia}(v)=s_{ia}(v)\ts s_{ia}(u)
\end{equation}
which implies $[s_{ia}^{(k)}, s_{ia}^{(r)}]=0$ for all $k,r$.
This completes the proof. \epf

\bth\label{thm:injaff}
The abstract algebra $\cal S$
generated by elements $s_{ij}^{(r)}$ with $i,j=1,\dots,N$ and $r\geq 0$
with the defining relations \eqref{sijoaff}--\eqref{rsrsaff}
is isomorphic to
$\Y^{\rm tw}_q(\oa_N)$.
\eth

\bpf We use the argument of the proof of Theorem~\ref{thm:inj}
appropriately modified for the affine case. Namely, we check first that the
matrix elements of the matrix $T(u)\ts \overline T(u^{-1})^t$ satisfy
the defining relations of $\cal S$ so that we have an algebra
homomorphism  ${\cal S}\to \Y^{\rm tw}_q(\oa_N)$. To prove its
injectivity regard $q$ as a formal variable and set
$V=\Y^{\rm tw}_q(\oa_N)$. Note
that the algebra $V_{\A}$ is generated by the elements
\beql{sigaff}
\sigma_{ij}^{(r)}=\frac{s_{ij}^{(r)}}{q-\qin},
\end{equation}
where the condition $i>j$ must hold if $r=0$.
The image of $\sigma_{ij}^{(r)}$ in $V_{\A}\ot_{\A}\C$
is $F_{ij}^{(r)}=E_{ij}\la^r-E_{ji}\la^{-r}$.
However, the elements
$F_{ij}^{(r)}$ span the fixed point subalgebra ${\wh\gl_N}^{\ts\theta}$
of $\wh\gl_N$
corresponding to the automorphism
\beql{thetaaff}
\theta: E_{ij}\ts \la^p\mapsto - E_{ji}\ts \la^{-p}.
\end{equation}
Thus, Lemma~\ref{lem:wpbwaff} and this observation complete the proof.
\epf

Consider the ordering on the generators of $\cal S$ defined in
Lemma~\ref{lem:wpbwaff}. We shall assume that
for the generators $s_{ij}^{(0)}$ the condition $i>j$ holds.

\bco\label{cor:pbwaff}
The ordered monomials in the generators
constitute a basis of the algebra $\cal S$.
\qed
\eco

\bco\label{cor:deform}
The subalgebra $\Y^{\rm tw}_q(\oa_N)$ of $\U_q(\wh\gl_N)$
specializes to the universal enveloping
algebra $\U(\wh\gl_N^{\ts\theta})$ as $q\to 1$.
\qed
\eco

The following proposition is immediate from \eqref{copraff}.

\bpr\label{prop:coidaff}
The subalgebra $\Y^{\rm tw}_q(\oa_N)$ of $\U_q(\wh\gl_N)$ is a left coideal
so that
\beql{Delsij}
\Delta(s_{ij}(u))=\sum_{k,l=1}^N t_{ik}(u)\ts\bar t_{jl}(u^{-1})\ot s_{kl}(u).
\end{equation}
\qed
\epr

The next proposition allows us to regard $\U^{\rm tw}_q(\oa_N)$ as
a subalgebra of $\Y^{\rm tw}_q(\oa_N)$.

\bpr\label{prop:embed} The map $s_{ij}\mapsto s_{ij}^{(0)}$
defines an embedding
$\U^{\rm tw}_q(\oa_N)\hookrightarrow\Y^{\rm tw}_q(\oa_N)$.
\epr

\bpf
It is clear from \eqref{sijoaff}--\eqref{rsrsaff}
that the elements $s_{ij}^{(0)}$ satisfy the defining relations for the $s_{ij}$;
see \eqref{sijo}--\eqref{rsrs}. Corollary~\ref{cor:pbwaff}
ensures that the homomorphism is injective. \epf

The following theorem establishes
the existence of the evaluation homomorphisms for the twisted
$q$-Yangians. We use notation \eqref{barS}. As before,
we combine the formal series $s_{ij}(u)$ into the matrix $S(u)$
so that $S(u)$ is a formal power series in $u^{-1}$
with matrix coefficients.

\bth\label{thm:eval} The mapping
\beql{evaluaff}
S(u)\mapsto S+\qin\ts u^{-1}\ts {\overline S}
\end{equation}
defines an algebra homomorphism $\Y^{\rm tw}_q(\oa_N)\to\U^{\rm tw}_q(\oa_N)$.
\eth

\bpf  Regarding $\Y^{\rm tw}_q(\oa_N)$ as an abstract algebra we need to verify
that the relation \eqref{rsrsaff} holds when we substitute
the right hand side of \eqref{evaluaff} for $S(u)$. That is, we need to verify
\begin{multline}\label{verirs}
(u\wt R-vR)(q\ts u\ts S_1+\overline S_1)(\wt R^t-uvR^t)
(q\ts v\ts S_2+\overline S_2)\\
=(q\ts v\ts S_2+\overline S_2)(\wt R^t-uvR^t)
(q\ts u\ts S_1+\overline S_1)(u\wt R-vR).
\end{multline}
Applying \eqref{defrel} and \eqref{relco} we derive the
relations
\beql{auxrsrs}
\bal
R\ts S_1 R^{\ts t} S_2=S_2R^{\ts t} S_1R,\qquad
\wt R\ts S_1 R^{\ts t} S_2=S_2R^{\ts t} S_1\wt R,\\
R\ts \overline S_1 R^{\ts t} S_2=S_2R^{\ts t}\overline S_1R,\qquad
\wt R\ts \overline S_1 \wt R^{\ts t} \overline S_2=
\overline S_2\wt R^{\ts t} \overline S_1\wt R,\\
R\ts \overline S_1 \wt R^{\ts t} S_2=S_2\wt R^{\ts t}\overline S_1R,\qquad
\wt R\ts S_1 \wt R^{\ts t}\overline S_2=\overline S_2\wt R^{\ts t} S_1\wt R,\\
R\ts \overline S_1\wt R^{\ts t} \overline S_2=
\overline S_2\wt R^{\ts t}\overline S_1R,\qquad
\wt R\ts S_1 R^{\ts t}\overline S_2=\overline S_2R^{\ts t} S_1\wt R.
\eal
\end{equation}
Expanding the products on the left and right hand sides of
\eqref{verirs} we see that in order
to complete the proof, it suffices to verify the following
four relations
\beql{fourrel}
\bal
R\ts S_1 R^{\ts t}\overline S_2-\wt R\ts\overline S_1 R^{\ts t} S_2
&=\overline S_2R^{\ts t} S_1R-S_2R^{\ts t}\overline S_1\wt R,\\
R\ts S_1 \wt R^{\ts t}\overline S_2-\wt R\ts\overline S_1 \wt R^{\ts t} S_2
&=\overline S_2\wt R^{\ts t} S_1R-S_2\wt R^{\ts t}\overline S_1\wt R,\\
\wt R\ts S_1 \wt R^{\ts t}S_2-q^{-2}\ts \wt R\ts\overline S_1 R^{\ts t}\overline S_2
&=S_2\wt R^{\ts t} S_1\wt R-q^{-2}\ts\overline S_2R^{\ts t}\overline S_1\wt R,\\
R\ts S_1\wt R^{\ts t} S_2-q^{-2}\ts R\ts\overline S_1 R^{\ts t} \overline S_2
&=S_2\wt R^{\ts t} S_1R-q^{-2}\ts \overline S_2R^{\ts t}\overline S_1 R.
\eal
\end{equation}
Let us prove the first of them.
Using \eqref{relRRt} replace $R$ and $\wt R$, respectively,
with $\wt R+(q-\qin)P$ and $R-(q-\qin)P$. Due to \eqref{auxrsrs},
we now have to check that
\beql{permrs}
P\ts S_1 R^{\ts t}\overline S_2+P\ts\overline S_1 R^{\ts t} S_2
=\overline S_2R^{\ts t} S_1P+S_2R^{\ts t}\overline S_1P.
\end{equation}
However, this is obvious if we observe that
$P\ts S_1=S_2P$, $P\ts \overline S_1=\overline S_2P$ and
$P\ts R^{\ts t}=R^{\ts t}P$. The proof of the second relation in \eqref{fourrel}
is the same. The third and forth relations are also verified in a way
similar to each other so we only consider the third one.
Note that by \eqref{relRRt},
\beql{qrrt}
R^{\ts t}=\wt R^{\ts t}+(q-\qin)\ts Q,\qquad
Q=P^{\ts t}=\sum_{i,j=1}^N E_{ij}\ot E_{ij}.
\end{equation}
Replacing $R^{\ts t}$ and $\wt R^{\ts t}$, respectively, with
$\wt R^{\ts t}+(q-\qin)\ts Q$ and $R^{\ts t}-(q-\qin)\ts Q$ in the third relation
and using again \eqref{auxrsrs} we conclude that it is now sufficient
to verify that
\beql{qrelaux}
\wt R\ts S_1 QS_2+q^{-2}\ts \wt R\ts\overline S_1 Q\overline S_2
=S_2 Q S_1\wt R+q^{-2}\ts\overline S_2Q\overline S_1\wt R.
\end{equation}
Since $PQ=QP=Q$ and
\beql{pmatp}
P\ts A_1=A_2 P
\end{equation}
for any matrix $A$,
we have
$P\ts\overline S_1 Q\overline S_2=\overline S_2Q\overline S_1P$.
Therefore, by \eqref{relRRt} we may replace \eqref{qrelaux} with
\beql{qrelaux2}
\wt R\ts S_1 QS_2+q^{-2}\ts R\ts\overline S_1 Q\overline S_2
=S_2 Q S_1\wt R+q^{-2}\ts\overline S_2Q\overline S_1 R.
\end{equation}
Finally, we have the chain of equalities
\beql{rsq}
\wt R\ts S_1 Q=\wt R\ts T_1 \overline T_1^{\ts t} Q =
\wt R\ts T_1 \overline T_2 Q =
\overline T_2 T_1 \wt R\ts Q    =\qin \overline T_2 T_1 Q =
\qin \overline T_2 T_2^{\ts t} Q =\qin \overline S_2 Q,
\end{equation}
where we have used \eqref{relco} and the relations
$\wt R\ts Q=Q \wt R=q^{-1}\ts Q$ together with the observation
that $A_1^t Q=A_2 Q$ implied by \eqref{pmatp}.  Thus,
$\wt R\ts S_1 QS_2=\qin \overline S_2 Q S_2$.
A similar argument shows that
\beql{rsqrfin}
R\ts\overline S_1 Q\overline S_2=q\ts S_2 Q \overline S_2,\qquad
S_2 Q S_1\wt R=\qin S_2 Q \overline S_2,\qquad
\overline S_2Q\overline S_1 R =q\ts \overline S_2Q S_2
\end{equation}
proving
\eqref{qrelaux2}
and the theorem. \epf

\subsection{Symplectic twisted $q$-Yangians}\label{subsec:sty}

As in Section~\ref{subsec:sc}, we use the matrix $G=(g_{ij})$ given by \eqref{g}.

\bde\label{def:symptwy}
The {\em twisted $q$-Yangian\/} $\Y^{\rm tw}_q(\spa_{2n})$
is defined as the subalgebra of $\U_q(\wh\gl_{2n})$ generated by the coefficients
$s_{ij}^{(r)}$ of
the matrix elements of the matrix $S(u)=T(u)\ts G\ts\overline T(u^{-1})^t$
and the elements $(s_{ii'}^{(0)})^{-1}$ with $i=1,3,\dots,2n-1$.
More precisely, the matrix elements are given by
\beql{sumatrs}
s_{ij}(u)=q\ts\sum_{a=1}^{n} t_{i,2a-1}(u)\ts \bar t_{j,2a}(u^{-1})
-\sum_{a=1}^{n} t_{i,2a}(u)\ts \bar t_{j,2a-1}(u^{-1})
\end{equation}
so that
\beql{scoeffs}
s_{ij}(u)=\sum_{r=0}^{\infty} s_{ij}^{(r)}\ts u^{-r}.
\end{equation}
\ede

Now we give a presentation of the algebra $\Y^{\rm tw}_q(\spa_{2n})$
in terms of generators and defining relations.
We shall use the notation \eqref{rtaff}.
Denote by $\cal S$ the associative algebra
with (abstract) generators
$s_{ij}^{(r)}$ with $1\leq i,j\leq 2n$ and $r\geq 0$,
$(s_{ii'}^{(0)})^{-1}$ with $i=1,3,\dots,2n-1$
and the following
defining relations written with the use of the generating series \eqref{scoeffs}
and the matrix $S(u)=(s_{ij}(u))$:
\beql{szeraffs}
s_{ij}^{(0)}=0\qquad\text{for}\quad i<j \quad\text{with}\quad   j\ne i';
\end{equation}
also, for any odd $i=1,3,\dots,2n-1$
\beql{qdss}
\bal
s_{i'i'}^{(0)}\ts s_{ii}^{(0)}-q^2\ts s_{i'i}^{(0)}\ts s_{ii'}^{(0)}=q^3,\\
s_{ii'}^{(0)}\ts (s_{ii'}^{(0)})^{-1}=(s_{ii'}^{(0)})^{-1}\ts s_{ii'}^{(0)}=1;
\eal
\end{equation}
and
\beql{rsrsaffs}
R(u,v)\ts S_1(u)\ts R^{\ts t}(u^{-1},v)\ts S_2(v)=S_2(v)
\ts R^{\ts t}(u^{-1},v)\ts S_1(u)\ts R(u,v).
\end{equation}

\bth\label{thm:injaffs}
The algebra $\cal S$
is isomorphic to the twisted $q$-Yangian
$\Y^{\rm tw}_q(\spa_{2n})$.
\eth

\bpf Let us verify first that the matrix $G$ satisfies
\beql{refaffg}
R(u,v)\ts G_1\ts R^{\ts t}(u^{-1},v)\ts G_2=G_2
\ts R^{\ts t}(u^{-1},v)\ts G_1\ts R(u,v).
\end{equation}
We have to check the four relations
\beql{relqsymp}
\bal
R\ts G_1 R^{\ts t} G_2&=G_2R^{\ts t} G_1R,\qquad
\wt R\ts G_1 R^{\ts t} G_2=G_2R^{\ts t} G_1\wt R,\\
R\ts G_1 \wt R^{\ts t} G_2&=G_2\wt R^{\ts t} G_1R,\qquad
\wt R\ts G_1 \wt R^{\ts t} G_2=
G_2\wt R^{\ts t} G_1\wt R.
\eal
\end{equation}
Applying \eqref{refg} together with \eqref{relRRt} and \eqref{qrrt}
we see that it is sufficient to check that
$R\ts G_1 Q G_2=G_2 Q G_1R$ which can be done by a direct calculation.

Further, repeating the arguments of the proof of Theorem~\ref{thm:injaff}
we construct the algebra homomorphism
${\cal S}\to\Y^{\rm tw}_q(\spa_{2n})$. Next we verify directly that
the elements
$s_{i'i'}^{(0)}\ts s_{ii}^{(0)}-q^2\ts s_{i'i}^{(0)}\ts s_{ii'}^{(0)}$
with odd $i$ are central in the algebra $\cal S$.

Now Lemma~\ref{lem:bass} and Lemma~\ref{lem:wpbwaff} imply that
the algebra $\cal S$ is spanned by the ordered monomials
in the generators $s_{ij}^{(r)}$. Here the triples $(r,i,j)$
are ordered lexicographically
except for the submonomials in the generators $s_{ij}^{(0)}$.
These submonomials should have the form \eqref{ordsymp}
with the $s_{ij}$ replaced by $s_{ij}^{(0)}$, respectively.

Finally, the linear independence of the ordered monomials is proved
in the same way as in Theorem~\ref{thm:injaff}. Namely,
setting $V=\Y^{\rm tw}_q(\spa_{2n})$ we consider
the corresponding elements of the algebra $V_{\A}$ given by
\beql{sigsyms}
\sigma_{ij}^{(0)}=\frac{s_{ij}^{(0)}-g_{ij}}{q-\qin},
\qquad\text{and}\qquad
\sigma_{ij}^{(r)}=\frac{s_{ij}^{(r)}}{q-\qin} \qquad\text{for}\quad r\geq 1,
\end{equation}
cf. \eqref{sigsym}.  Then we have modulo $(q-1)$
\beql{sigmatauss}
\sigma_{ij}^{(r)}\equiv \wt g_{j'j}\ts\tau_{ij'}^{(r)}+\wt g_{ii'}
\ts\bar\tau_{ji'}^{(r)},
\end{equation}
where
$\wt g_{ii'}=(-1)^{i-1}$; see \eqref{sigmataus}.
Hence, the image of $\sigma_{ij}^{(r)}$ in $V_{\A}\ot_{\A}\C$ is
\beql{fijrs}
F_{ij}^{(r)}:=\wt g_{j'j}\ts E_{ij'}\la^r-\wt g_{ii'}\ts E_{ji'}\la^{-r}.
\end{equation}
The elements $F_{ij}^{(r)}$ span a fixed point
subalgebra ${\wh\gl_{2n}}^{\ts\theta}$
of $\wh\gl_{2n}$ with respect to the automorphism
\beql{autosyms}
\theta: E_{ij}\la^r\mapsto (-1)^{i+j-1} E_{j'i'}\la^{-r}.
\end{equation}
The application of the Poincar\'e--Birkhoff--Witt theorem
to the Lie algebra ${\wh\gl_{2n}}^{\ts\theta}$ completes the proof.
\epf

Consider the ordering on the generators of $\cal S$ defined in
the proof of Theorem~\ref{thm:injaffs}.

\bco\label{cor:pbwaffs}
The ordered monomials in the generators
constitute a basis of the algebra $\cal S$.
\qed
\eco

\bco\label{cor:deforms}
The subalgebra $\Y^{\rm tw}_q(\spa_{2n})$ of $\U_q(\wh\gl_{2n})$
specializes to the universal enveloping algebra $\U(\wh\gl_{2n}^{\ts\theta})$
as $q\to 1$.
\qed
\eco

The following proposition is immediate from \eqref{copraff}.

\bpr\label{prop:coidaffs}
The subalgebra $\Y^{\rm tw}_q(\spa_{2n})$ of $\U_q(\wh\gl_{2n})$ is a left coideal
so that
\beql{Delsijs}
\Delta(s_{ij}(u))=\sum_{k,l=1}^{2n} t_{ik}(u)\ts\bar t_{jl}(u^{-1})\ot s_{kl}(u).
\end{equation}
\qed
\epr

The next proposition allows us to regard $\U^{\rm tw}_q(\spa_{2n})$ as
a subalgebra of $\Y^{\rm tw}_q(\spa_{2n})$.

\bpr\label{prop:embeds} The map
\beql{embedss}
s_{ij}\mapsto s_{ij}^{(0)},\qquad s_{ii'}^{-1}\mapsto (s_{ii'}^{(0)})^{-1}
\end{equation}
defines an embedding
$\U^{\rm tw}_q(\spa_{2n})\hookrightarrow\Y^{\rm tw}_q(\spa_{2n})$.
\epr

\bpf
By the defining relations in $\U^{\rm tw}_q(\spa_{2n})$ and
$\Y^{\rm tw}_q(\spa_{2n})$ the mapping is clearly an algebra homomorphism.
Corollary~\ref{cor:pbwaffs}
ensures that the homomorphism is injective. \epf

The following is an analog of Theorem~\ref{thm:eval}
for the symplectic case. We use notation \eqref{barSs}.

\bth\label{thm:evals} The mapping
\beql{evaluaffs}
S(u)\mapsto S+q\ts u^{-1}\ts {\overline S},\qquad
(s_{ii'}^{(0)})^{-1}\mapsto  s_{ii'}^{-1}
\end{equation}
defines an algebra homomorphism $\Y^{\rm tw}_q(\spa_{2n})\to\U^{\rm tw}_q(\spa_{2n})$.
\eth

\bpf By analogy with the orthogonal case (see Theorem~\ref{thm:eval}),
we need to verify
\begin{multline}\label{verirss}
(u\wt R-vR)(u\ts S_1+q\ts \overline S_1)(\wt R^t-uvR^t)
(v\ts S_2+q\ts \overline S_2)\\
=(v\ts S_2+q\ts \overline S_2)(\wt R^t-uvR^t)
(u\ts S_1+q\ts \overline S_1)(u\wt R-vR).
\end{multline}
Note that \eqref{auxrsrs} and the first two relations in
\eqref{fourrel} still hold in the same form for the symplectic case.
To prove the theorem we shall verify that the third and fourth
relations in \eqref{fourrel} are respectively replaced by
\beql{fourrels}
\bal
\wt R\ts S_1 \wt R^{\ts t}S_2-q^{2}\ts \wt R\ts\overline S_1 R^{\ts t}\overline S_2
&=S_2\wt R^{\ts t} S_1\wt R-q^{2}\ts\overline S_2R^{\ts t}\overline S_1\wt R,\\
R\ts S_1\wt R^{\ts t} S_2-q^{2}\ts R\ts\overline S_1 R^{\ts t} \overline S_2
&=S_2\wt R^{\ts t} S_1R-q^{2}\ts \overline S_2R^{\ts t}\overline S_1 R.
\eal
\end{equation}
The chain of equalities \eqref{rsq} is replaced with
\begin{multline}\label{rsqs}
\wt R\ts S_1 Q=\wt R\ts T_1 G_1\overline T_1^{\ts t} Q =
\wt R\ts T_1 \overline T_2 G_1 Q =
\overline T_2 T_1 \wt R G_1\ts Q\\
=-q\ts \overline T_2 T_1 G_2 Q =
-q\ts \overline T_2 G_2 T_2^{\ts t} Q =-q\ts \overline S_2 Q,
\end{multline}
where we have used \eqref{relco} and the relation
$\wt R G_1\ts Q=-q\ts\ts G_2 Q$ which is verified directly.
Thus,
$\wt R\ts S_1 QS_2=-q\ts \overline S_2 Q S_2$.
A similar argument shows that
\beql{rsqrfins}
R\ts\overline S_1 Q\overline S_2=-\qin\ts S_2 Q \overline S_2,\qquad
S_2 Q S_1\wt R=-q\ts S_2 Q \overline S_2,\qquad
\overline S_2Q\overline S_1 R =-\qin\ts \overline S_2Q S_2
\end{equation}
completing the proof. \epf

\subsection{Comments on possible variations of the definitions}\label{subsec:com}

$\bullet$\quad By analogy with the algebras $\U^{\rm tw}_q(\oa_{N})$
and $\U^{\rm tw}_q(\spa_{2n})$
we can introduce the matrices   $\overline S(u) =\big(\bar s_{ij}(u)\big)$ by
\beql{barSuorth}
\overline S(u) =\overline  T(u)T(u^{-1})^t\qquad\text{and}\qquad
\overline S(u) =\overline  T(u)\ts G\ts T(u^{-1})^t
\end{equation}
in the orthogonal and symplectic case, respectively; cf. Definitions~\ref{def:orttwy}
and \ref{def:symptwy}.
Then one can derive the following relations form \eqref{defrelaff}:
\begin{multline}\label{relsbarsorth}
(uq-u^{-1}q^{-1})\ts\bar s_{ij}(u) = \\
(uq^{\delta_{ij}}-u^{-1}q^{-\delta_{ij}})\, s_{ji}(u^{-1})
+(q-q^{-1})(u\delta_{j<i}+u^{-1}\delta_{i<j})\ts s_{ij}(u^{-1})
\end{multline}
in the orthogonal case, and
\begin{multline}\label{relsbarssympl}
(u^{-1}q-uq^{-1})\ts\bar s_{ij}(u) = \\
(uq^{\delta_{ij}}-u^{-1}q^{-\delta_{ij}})\, s_{ji}(u^{-1})
+(q-q^{-1})(u\delta_{i<j}+u^{-1}\delta_{j<i})\ts s_{ij}(u^{-1})
\end{multline}
in the symplectic case. This implies that
the coefficients of the series
$\bar s_{ij}(u)$ generate the same subalgebra
$\Y^{\rm tw}_q(\oa_{N})$ or $\Y^{\rm tw}_q(\spa_{2n})$ of $\U_q(\wh\gl_N)$.

\bigskip
\noindent
$\bullet$\quad The construction of the coideal subalgebras
$\Y^{\rm tw}_q(\oa_{N})$ and $\Y^{\rm tw}_q(\spa_{2n})$
can be generalized to the case of the
centrally extended quantum affine algebra
$\U_{q}(\wh{\gl}_{N})_{c}$. We use its presentation given in
\cite{d:sr}. The $R$-matrix $R(u)$ we employ here is related to $R(u,v)$
by
\beql{ruruv}
(uq^{-1}-q)\,R(u)=R(u,1).
\end{equation}
The defining relations \eqref{defrelaff} are then replaced by
\beql{defrelcent}
\bal
R(u/v)\ts T_{1}(u)\ts T_{2}(v) &=
T_{2}(v)\ts T_{1}(u)\ts R(u/v),\\
R(u/v)\ts \overline T_{1}(u)\ts \overline T_{2}(v) &=
\overline T_{2}(v)\ts \overline T_{1}(u)\ts R(u/v),\\
R(u\tss q^{-c}/v)\ts \overline T_{1}(u)\ts T_{2}(v) &=
T_{2}(v)\ts \overline T_{1}(u)\ts R(u\tss q^{c}/v).
\eal
\end{equation}
The definition of the matrix $S(u)$ is modified to
\beql{sucntrext}
S(u)=T(uq^{-c})\ts\overline T(u^{-1})^t\qquad\text{and}\qquad
S(u)=T(uq^{-c})\ts G\ts \overline T(u^{-1})^t
\end{equation}
in the orthogonal and symplectic case, respectively.
However, one can demonstrate that $S(u)$ still satisfies
the same reflection equation \eqref{rsrsaff} or \eqref{rsrsaffs}
which does not involve the central charge $c$.
The corresponding subalgebra is therefore isomorphic
to the twisted $q$-Yangian, as an abstract algebra.

\section{Centers of the twisted $q$-Yangians}\label{sec:center}
\setcounter{equation}{0}

In this section we construct a formal series whose coefficients
belong to the center of the twisted $q$-Yangian.

We start by recalling the well-known construction of the quantum determinants
for the quantum affine algebra $\U_q(\wh\gl_N)$; see e.g.
\cite{c:ni}, \cite{j:qu}, \cite{rtf:ql}.
We use an approach analogous to the case of the
Yangian $\Y(\gl_N)$ \cite{ik:lm}, \cite{ks:qs}; see also
\cite{mno:yc} for a detailed exposition.

Let us consider
the multiple tensor product
$\U_q(\wh\gl_N)\ot (\End\C^N)^{\ot\tss r}$ and use the notation of \eqref{defrelaff}.
Then we have the following corollary of \eqref{YBE} and \eqref{defrelaff}:
\beql{fundam}
R(u_1,\dots,u_r)\ts T_1(u_1)\cdots T_r(u_r)=T_r(u_r)
\cdots T_1(u_1)\ts  R(u_1,\dots,u_r),
\end{equation}
where
\beql{Rlong}
R(u_1,\dots,u_r)=\prod_{i<j}R_{ij}(u_i,u_j),
\end{equation}
with the product taken in the lexicographical order on the pairs $(i,j)$.
The proof of \eqref{fundam} is exactly the
same as for the Yangians; see e.g. \cite{mno:yc}.
Furthermore, consider the $q$-permutation operator
$P^{\tss q}\in\End(\C^N\ot\C^N)$
defined by
\beql{qperm}
P^{\tss q}=\sum_{i}E_{ii}\ot E_{ii}+ q\tss\sum_{i> j}E_{ij}\ot
E_{ji}+ \qin\sum_{i< j}E_{ij}\ot E_{ji}.
\end{equation}
The action of symmetric group $\Sym_r$ on the space $(\C^N)^{\ot\tss r}$
can be defined by setting $s_i\mapsto P^{\tss q}_{s_i}:=
P^{\tss q}_{i,i+1}$ for $i=1,\dots,r-1$,
where $s_i$ denotes the transposition $(i,i+1)$.
If $\sigma=s_{i_1}\cdots s_{i_l}$ is a reduced decomposition
of an element $\sigma\in \Sym_r$ we set
$P^{\tss q}_{\sigma}=P^{\tss q}_{s_{i_1}}\cdots P^{\tss q}_{s_{i_l}}$.
We denote by $A^q_r$ the $q$-antisymmetrizer
\beql{antisym}
A^q_r=\sum_{\sigma\in\Sym_r}\sgn\ts\sigma\cdot P^{\tss q}_{\sigma}.
\end{equation}

The following proposition is proved by induction on $r$ in the same way as for the
Yangians \cite{mno:yc} with the use of a property of the reduced decompositions
\cite[p.\tss50]{h:il}.

\bpr\label{prop:antisym}
We have the relation in $\ts\End (\C^N)^{\ot\tss r}$:
\beql{ranti}
R(1,q^{-2},\dots,q^{-2r+2})=\prod_{0\leq i<j\leq r-1}(q^{-2i}-q^{-2j})\ts A^q_r.
\end{equation}
\epr

Now \eqref{fundam} implies
\beql{anitt}
A^q_r\ts T_1(u)\cdots T_r(q^{-2r+2}u)=T_r(q^{-2r+2}u)\cdots T_1(u)\ts  A^q_r
\end{equation}
which equals
\beql{matelmi}
\sum_{a_i,b_i}{t\ts}^{a_1\cdots\ts a_r}_{b_1\cdots\ts b_r}(u)\ot
E_{a_1b_1}\ot\cdots\ot E_{a_rb_r}
\end{equation}
for some elements ${t\ts}^{a_1\cdots\ts a_r}_{b_1\cdots\ts b_r}(u)$
of $\U_q(\wh\gl_N)[[u^{-1}]]$
which we call the {\it quantum minors\/}.
They can be given by the following formulas which are immediate from the definition.
If $a_1<\cdots<a_r$ then
\beql{qminorgen}
{t\ts}^{a_1\cdots\ts a_r}_{b_1\cdots\ts b_r}(u)=
\sum_{\sigma\in \Sym_r} (-q)^{-l(\sigma)} \cdot t_{a_{\sigma(1)}b_1}(u)\cdots
t_{a_{\sigma(r)}b_r}(q^{-2r+2}u),
\end{equation}
and for any $\tau\in\Sym_r$ we have
\beql{qmsym}
{t\ts}^{a_{\tau(1)}\cdots\ts a_{\tau(r)}}_{b_1\cdots\ts b_r}(u)=
(-q)^{l(\tau)}{t\ts}^{a_1\cdots\ts a_r}_{b_1\cdots\ts b_r}(u),
\end{equation}
where $l(\sigma)$ denotes the length of the permutation $\sigma$.
If $b_1<\cdots<b_r$ (and the $a_i$ are arbitrary) then
\beql{qminorgen2}
{t\ts}^{a_1\cdots\ts a_r}_{b_1\cdots\ts b_r}(u)=
\sum_{\sigma\in \Sym_r} (-q)^{l(\sigma)} \cdot t_{a_rb_{\sigma(r)}}(q^{-2r+2}u)\cdots
t_{a_1b_{\sigma(1)}}(u),
\end{equation}
and for any $\tau\in\Sym_r$ we have
\beql{qmsym2}
{t\ts}^{a_1\cdots\ts a_r}_{b_{\tau(1)}\cdots\ts b_{\tau(r)}}(u)=
(-q)^{-l(\tau)}{t\ts}^{a_1\cdots\ts a_r}_{b_1\cdots\ts b_r}(u).
\end{equation}
Moreover, the quantum minor is zero if two top or two bottom indices
are equal.
Note also that the standard row and column expansion
formulas can be easily derived from \eqref{qminorgen} or
\eqref{qminorgen2}. In particular, we have
\beql{columnexp}
{t\ts}^{a_1\cdots\ts a_{r}}_{b_1\cdots\ts b_{r}}(u)
=\sum_{l=1}^{r}(-q)^{-l+1}\ts{t\ts}^{}_{a_lb_1}(u)\ts
{t\ts}^{a_1\cdots\wh{a}_l\cdots\ts a_{r}}_{b_2\cdots\ts b_{r}}(q^{-2}\ts u),
\end{equation}
where the hat indicates the index to be omitted.

The quantum minors
${\bar t\ts}^{a_1\cdots\ts a_r}_{b_1\cdots\ts b_r}(u)$
of the matrix $\overline T(u)$ are given by the same formulas
where the $t_{ij}(u)$ are respectively replaced with $\bar t_{ij}(u)$.
Furthermore, for any indices $i,j$ we have the well known
relations which are deduced from \eqref{fundam}:
\beql{center}
[{t}^{}_{c_id_j}(u),
{t\ts}^{c_1\cdots\ts c_r}_{d_1\cdots\ts d_r}(v)]=0,\qquad
[{t}^{}_{c_id_j}(u),
{\bar t\ts}^{c_1\cdots\ts c_r}_{d_1\cdots\ts d_r}(v)]=0,
\end{equation}
and the same holds with ${t}^{}_{c_id_j}(u)$ replaced by
${\bar t}^{}_{c_id_j}(u)$.
For their proof introduce an extra copy of
$\End\C^N$ as a tensor factor which will be enumerated by the index $0$.
Now we specialize
the parameters $u_i$ in \eqref{fundam} as follows:
\beql{specpar}
u_0=v, \qquad u_i=q^{-2i+2}u\quad\text{for}\ \ i=1,\dots,r.
\end{equation}
Then by Proposition~\ref{prop:antisym}
the element \eqref{Rlong} will take the form
\beql{RA}
R(v,u,\dots,q^{-2r+2}u)=\prod_{i=1}^r R_{0i}(v,q^{-2i+2}u)\ts A^q_r.
\end{equation}
Using the definition of the quantum minors
and equating the matrix elements on both sides of \eqref{fundam}
we get the first relation in \eqref{center}; cf. \cite{mno:yc}.
The proof of the second is similar.

The {\it quantum determinants\/}
of the matrices $T(u)$ and $\overline T(u)$ are respectively
defined by the relations
\beql{qdets}
\qdet T(u)={t\ts}^{1\cdots\ts N}_{1\cdots\ts N}(u),\qquad
\qdet \overline T(u)={\bar t\ts}^{1\cdots\ts N}_{1\cdots\ts N}(u).
\end{equation}
Write
\beql{coeffqdet}
\qdet T(u)=\sum_{k=0}^{\infty} d_k\ts u^{-k},\qquad
\qdet \overline T(u)=\sum_{k=0}^{\infty} \bar d_k\ts u^{k},
\qquad d_k,\bar d_k\in \U_q(\wh\gl_N).
\end{equation}
Due to the property \eqref{center}, the coefficients $d_k$ and $\bar d_k$
belong to the center of
the algebra $\U_q(\wh\gl_N)$.
Note also the relation $d_0\bar d_0=1$
which is implied by \eqref{defrelaff}.

\subsection{Sklyanin determinant}\label{subsec:sdet}

We shall consider the orthogonal and symplectic cases
simultaneously unless otherwise stated. The symbol $\Y^{\rm tw}_q$
will denote either $\Y^{\rm tw}_q(\oa_{N})$ or $\Y^{\rm tw}_q(\spa_{N})$
(the latter with $N=2n$).
As before, we denote by $S(u)$ the matrix $T(u)G\ts \overline T(u^{-1})^t$,
where $G$ is the identity matrix in the orthogonal case, and
$G$ is given by \eqref{g} in the symplectic case.
Then $S(u)$ satisfies the relations
\eqref{rsrsaff} and \eqref{rsrsaffs}, respectively, which have the same form.
Our arguments are similar to those used in \cite{o:ty} and \cite{mno:yc}
for the construction of the Sklyanin determinant for the twisted Yangians.

The following relation in the algebra $\Y^{\rm tw}_q\ot (\End\C^N)^{\ot\tss r}$
is a corollary of \eqref{rsrsaffs}.
\begin{multline}\label{fundamtw}
R(u_1,\dots,u_r)\ts S_1(u_1)R_{12}^{\tss t}\cdots R_{1r}^{\tss t}   S_2(u_2)
R_{23}^{\tss t}\cdots R_{2r}^{\tss t}   S_3(u_3) \cdots R_{r-1,r}^{\tss t}
S_r(u_r)=\\
S_r(u_r)R_{r-1,r}^{\tss t}\cdots S_3(u_3)   R_{2r}^{\tss t} \cdots R_{23}^{\tss t}
S_2(u_2)R_{1r}^{\tss t} \cdots R_{12}^{\tss t} S_1(u_1)\ts R(u_1,\dots,u_r),
\end{multline}
where $R(u_1,\dots,u_r)$ is given by
\eqref{Rlong} and $R_{ij}^{\tss t}=R_{ij}^{\tss t}(u_i^{-1},u_j)$.
Now take $r=N$ and specify $u_i=u\ts q^{-2i+2}$.
Then Proposition~\ref{prop:antisym} implies
\begin{multline}\label{sdetmatr}
A^q_N\ts S_1(u)R_{12}^{\tss t}\cdots R_{1N}^{\tss t}    S_2(u\ts q^{-2})
R_{23}^{\tss t}\cdots R_{2N}^{\tss t}   S_3(u\ts q^{-4}) \cdots R_{N-1,N}^{\tss t}
S_N(u\ts q^{-2N+2})=\\
S_N(u\ts q^{-2N+2})R_{N-1,N}^{\tss t}\cdots S_3(u\ts q^{-4})    R_{2N}^{\tss t}
\cdots R_{23}^{\tss t}
S_2(u\ts q^{-2})R_{1N}^{\tss t} \cdots R_{12}^{\tss t} S_1(u)\ts A^q_N.
\end{multline}
Since the $q$-antisymmetrizer $A^q_N$ is proportional to an idempotent
and maps the space $(\C^N)^{\ot\tss N}$ into
a one-dimensional subspace,
both sides must be equal to $A^q_N$ times a series $\sdet S(u)$ in $u^{-1}$
with coefficients in
$\Y^{\rm tw}_q$. We call this series the {\it Sklyanin determinant\/}.

The following theorem provides an expression of $\sdet S(u)$ in terms
of the quantum determinants.

\bth\label{thm:sdetqdet}
We have
\beql{sdetqdet}
\sdet S(u)=\ga_N(u)\ts\qdet T(u)\ts\qdet \overline T(q^{2N-2}u^{-1}),
\end{equation}
where
\beql{gau}
\ga_N(u)=\prod_{1\leq i<j \leq N} (q^{2i-2}u^{-1}-q^{-2j+2}u)\cdot
\begin{cases}
1 &\quad\text{in the case of}\quad \oa_N,\\
\dfrac{q^{n-2}-q^n u^2}{q^{2n-2}-q^{-2n} u^2}
&\quad\text{in the case of}\quad \spa_{2n}.
\end{cases}
\end{equation}
\eth

\bpf We follow the arguments of \cite[Section~4]{mno:yc}.
Substitute $S(u)=T(u)G\ts \overline T(u^{-1})^t$ into \eqref{sdetmatr}
and transform the left hand side using the relations
\beql{rrtrt}
\overline T_i(u^{-1})^t\ts R_{ij}^t(u^{-1},v)\ts T_j(v)=T_j(v)
\ts R_{ij}^t(u^{-1},v)\ts \overline T_i(u^{-1})^t
\end{equation}
which are implied by \eqref{defrelaff}. We then bring it to the form
\beql{attt}
A^q_N \ts T_1(u)\cdots T_N(q^{-2N+2}u)\ts \wt R(u)\ts
\overline T_1(u^{-1})^t\cdots \overline T_N(q^{2N-2}u^{-1})^t,
\end{equation}
where
\beql{rtuu}
\wt R(u)=G_1R_{12}^{\tss t}\cdots R_{1N}^{\tss t}   G_2
R_{23}^{\tss t}\cdots R_{2N}^{\tss t}   G_3 \cdots R_{N-1,N}^{\tss t} G_N.
\end{equation}
By the definition of the quantum determinant we have
\beql{atttqdet}
A^q_N \ts T_1(u)\cdots T_N(q^{-2N+2}u)=A^q_N \ts \qdet T(u).
\end{equation}
Further, using the homomorphism $S(u)\mapsto G$
we derive from \eqref{sdetmatr}
\beql{agrgr}
A^q_N \wt R(u) =
G_NR_{N-1,N}^{\tss t} \cdots
G_3 R_{2N}^{\tss t}\cdots R_{23}^{\tss t}
G_2 R_{1N}^{\tss t}\cdots   R_{12}^{\tss t}
G_1 A^q_N.
\end{equation}
Thus, this expression equals $A^q_N \ga_N(u)$ for a scalar function $\ga_N(u)$.
Using the explicit formulas \eqref{qminorgen} one easily derives that
\beql{atttqdett}
A^q_N \ts \overline T_1(u^{-1})^t\cdots \overline T_N(q^{2N-2}u^{-1})^t
=A^q_N \ts \qdet \overline T(q^{2N-2}u^{-1}).
\end{equation}
It remains to calculate the function $\ga_N(u)$.
Consider first the orthogonal case. Apply the operator $A^q_N \wt R(u)$
to the basis vector $e_1\ot\cdots\ot\ts e_N$, where the $e_i$ denote the
canonical basis vectors of $\C^N$. By \eqref{rtaff} the result
is clearly $\ga_N(u)\ts A^q_N (e_1\ot\cdots\ot e_N)$
with $\ga_N(u)$ given by \eqref{gau}. In the symplectic case apply
$A^q_{2n} \wt R(u)$ to the basis vector
\beql{basvects}
v=e_{2n-1}\ot e_{2n-3}\ot\cdots\ot e_1\ot e_2\ot e_4
\ot\cdots\ot e_{2n}.
\end{equation}
Since $G\ts e_{2k}=q\ts e_{2k-1}$ the vector $A^q_{2n} \wt R(u)\ts v$
equals
\beql{halfant}
q^n\ts \prod_{n\leq i<j \leq 2n} (q^{2i-2}u^{-1}-q^{-2j+2}u)\ts
A^q_{2n} G_1R_{12}^{\tss t}\cdots R_{1,2n}^{\tss t}     G_2\cdots G_n
R_{n,n+1}^{\tss t}\cdots R_{n,2n}^{\tss t}\ts w,
\end{equation}
where
\beql{vectw}
w=e_{2n-1}\ot e_{2n-3}\ot\cdots\ot e_1\ot e_1\ot e_3
\ot\cdots\ot e_{2n-1}.
\end{equation}
If $j>n+1$ then $R_{n,j}^{\tss t}\ts w=(q^{2n-2}u^{-1}-q^{-2j+2}u)\ts w$.
Further, we have
\beql{ree}
R_{n,n+1}^{\tss t}\ts (e_1\ot e_1)=
(q^{2n-3}u^{-1}-q^{-2n+1}u)\ts (e_1\ot e_1)+(\qin-q)\ts q^{-2n}\ts u
\sum_{k=2}^{2n}(e_k\ot e_k).
\end{equation}
Next apply the operator $G_n$ and note that
due to the subsequent application of the $q$-antisymmetrizer, we may
only keep the linear combination of the tensor products containing
$e_1$ or $e_2$ on the $n$-th and $(n+1)$-th places (see \cite[Section~4]{mno:yc}
for a similar argument in the symplectic twisted Yangian case).
That is, we may write
\begin{multline}\label{fsteparw}
A^q_{2n} G_1R_{12}^{\tss t}\cdots R_{1,2n}^{\tss t}     G_2\cdots G_n
R_{n,n+1}^{\tss t}\ts w=\\
(q^{2n-4}u^{-1}-q^{-2n+2}u)\ts
A^q_{2n} G_1R_{12}^{\tss t}\cdots R_{1,2n}^{\tss t}     G_2\cdots G_{n-1}
R_{n-1,n}^{\tss t}\cdots R_{n-1,2n}^{\tss t}\ts w'
\end{multline}
with
\beql{vectwpr}
w'=e_{2n-1}\ot e_{2n-3}\ot\cdots\ot e_3\ot e_1\ot e_2\ot e_3
\ot\cdots\ot e_{2n-1}.
\end{equation}
By the skew-symmetry of the antisymmetrizer replace $w'$ with
the vector
\beql{vectwprpr}
w''=e_{2n-1}\ot e_{2n-3}\ot\cdots\ot e_3\ot e_3\ot e_2\ot e_1
\ot\cdots\ot e_{2n-1},
\end{equation}
taking the sign into account. Continuing the calculation in a similar manner,
we conclude that the product of the scalar factors occurring in the procedure
will coincide with $\ga_{2n}(u)$ given in \eqref{gau} with $N=2n$.
\epf

The centrality of $\qdet T(u)$ and $\qdet \overline T(u)$ in $\U_q(\wh\gl_N)$
immediately implies the corresponding property of the
Sklyanin determinant $\sdet S(u)$.

\bco\label{cor:centersdet}
The coefficients of the series $\sdet S(u)$ belong to the center
of the algebra $\Y^{\rm tw}_q$.
\qed
\eco

Introduce the series $c(u)$
and the elements $c_k$ of the center of the algebra $\Y^{\rm tw}_q$
by the formula
\beql{expansdet}
c(u)=\ga_N(u)^{-1}\ts\sdet S(u)=1+\sum_{k=1}^{\infty} c_k\ts u^{-k}.
\end{equation}
The series begins with $1$ due to Theorem~\ref{thm:sdetqdet}
and the relation $d_0\bar d_0=1$.

\bpr\label{prop:algindep}
The coefficients $c_k$, $k\geq 1$ are algebraically independent.
\epr

\bpf The coefficients $\{d_k,\  k\geq 0,\quad \bar d_k,\  k\geq 1\}$
of the quantum determinants are algebraically independent.
This can be derived by analogy with the case of the Yangian;
see e.g. \cite[Section~2]{mno:yc}.   The key observation here
is the isomorphism \eqref{limim}.
The statement is now implied by Theorem~\ref{thm:sdetqdet}.
\epf

Since both evaluation homomorphisms \eqref{evaluaff} and \eqref{evaluaffs}
are surjective, we obtain families of central elements in the algebras
$\U^{\rm tw}_q(\oa_{N})$ and $\U^{\rm tw}_q(\spa_{2n})$
as images of the coefficients of $\sdet S(u)$. In other words, we have
the following result.

\bpr\label{prop:sklimages}
The coefficients of the Sklyanin determinants
$\sdet (S+q^{-1}\tss u^{-1}\tss\overline S)$ and
$\sdet (S+q\tss u^{-1}\tss\overline S)$ are central elements in the algebras
$\U^{\rm tw}_q(\oa_{N})$ and $\U^{\rm tw}_q(\spa_{2n})$,
respectively.
\epr

In what follows we only consider the case of the orthogonal
twisted $q$-Yangian $\Y^{\rm tw}_q(\oa_{N})$.
In order to produce an explicit formula for the corresponding
Sklyanin determinant we introduce a map
\beql{mapp}
\pi^{}_N:\Sym_{N}\to\Sym_N,\qquad p\mapsto p'
\end{equation}
which was previously used in the `short'
formula for the Sklyanin determinant for the twisted
Yangians; see \cite{m:ya}. This map
is defined by an inductive procedure.
Given a set of positive integers
$\om_1<\cdots<\om_N$ we
regard $\Sym_N$ as the group of their permutations.
If $N=2$ we define $\pi^{}_2$ as the map $\Sym_2\to \Sym_{2}$
whose image is the identity permutation.
For $N>2$ define a map from the set of ordered pairs $(\om_k,\om_l)$
with $k\ne l$ into itself by the rule
\beql{ordpair}
\begin{alignedat}{2}
(\om_k,\om_l)&\mapsto (\om_l,\om_k),&&\qquad k,l<N,\\
(\om_k,\om_N)&\mapsto (\om_{N-1},\om_k),&&\qquad k<N-1,\\
(\om_N,\om_k)&\mapsto (\om_k,\om_{N-1}),&&\qquad k<N-1,\\
(\om_{N-1},\om_N)&\mapsto (\om_{N-1},\om_{N-2}),\\
(\om_{N},\om_{N-1})&\mapsto (\om_{N-1},\om_{N-2}).
\end{alignedat}
\end{equation}
Let $p=(p^{}_1,\dots,p^{}_N)$ be a permutation of the indices
$\om_1,\dots,\om_N$. Its image under
the map $\pi^{}_N$
is the permutation $p_{}^{\ts\prime}=
(p^{\ts\prime}_1,\dots,p^{\ts\prime}_{N-1},\om_N)$, where the pair
$(p^{\ts\prime}_1,p^{\ts\prime}_{N-1})$ is the image
of the ordered pair $(p^{}_1,p^{}_N)$ under the map \eqref{ordpair}.
Then the pair $(p^{\ts\prime}_2,p^{\ts\prime}_{N-2})$ is found as
the image of $(p^{}_2,p^{}_{N-1})$
under the map \eqref{ordpair} which is defined on the set
of ordered pairs of elements obtained from $(\om_1,\dots,\om_N)$
by deleting $p^{}_1$ and $p^{}_N$; etc.
The map $\pi^{}_N$
has curious combinatorial properties which were observed
by Lascoux; see \cite{m:sp}.
In particular,
each fiber of this map is an interval in $\Sym_N$
with respect to the Bruhat order, isomorphic to a Boolean poset.

\bex
The Bruhat order on $\Sym_3$ and the fibers of the map $\pi_3$:

\bigskip

\begin{center}
\begin{picture}(400,100)
\thinlines

\put(50,0){\circle*{3}}
\put(20,30){\circle*{3}}
\put(80,30){\circle*{3}}
\put(20,60){\circle*{3}}
\put(80,60){\circle*{3}}
\put(50,90){\circle*{3}}

\put(50,0){\line(1,1){30}}
\put(50,0){\line(-1,1){30}}
\put(20,60){\line(1,1){30}}
\put(80,60){\line(-1,1){30}}
\put(20,30){\line(0,1){30}}
\put(80,30){\line(0,1){30}}

\put(20,30){\line(2,1){60}}
\put(80,30){\line(-2,1){60}}

\put(45,-10){\scriptsize $123$ }
\put(0,25){\scriptsize $213$ }
\put(85,25){\scriptsize $132$ }
\put(0,60){\scriptsize $231$ }
\put(85,60){\scriptsize $312$ }
\put(45,95){\scriptsize $321$ }

\put(180,20){\circle*{3}}
\put(180,70){\circle*{3}}
\put(155,45){\circle*{3}}
\put(205,45){\circle*{3}}

\put(180,20){\line(1,1){25}}
\put(180,20){\line(-1,1){25}}
\put(155,45){\line(1,1){25}}
\put(205,45){\line(-1,1){25}}

\put(175,10){\scriptsize $123$ }
\put(135,45){\scriptsize $213$ }
\put(175,75){\scriptsize $312$ }
\put(210,45){\scriptsize $132$ }

\put(230,45){$\longrightarrow\ \  $\scriptsize$213$ }

\put(330,30){\circle*{3}}
\put(330,60){\circle*{3}}

\put(330,30){\line(0,1){30}}

\put(325,15){\scriptsize $231$ }
\put(325,70){\scriptsize $321$ }

\put(350,45){$\longrightarrow\ \ $\scriptsize$ 123$ }

\end{picture}
\end{center}
\eex

\vskip0.5cm

Consider the matrix
$\overline S(u)$ introduced in \eqref{barSuorth}.
By \eqref{relsbarsorth} the matrix elements of $\overline S(u)$
are formal series in $u$ with coefficients in
the subalgebra $\Y^{\rm tw}_q(\oa_{N})$.

Let $n$ denote the rank of the Lie algebra $\oa_N$
so that $N=2n$ or $N=2n+1$.

\bth\label{thm:exp} We have an explicit formula
\beql{ef}
\bal
c(u)=\sum_{p\in \Sym_{N}} (-q)^{-l(p)+l(p')}
\ts\bar s^{\ts t}_{p^{}_1p^{\ts\prime}_1}(u^{-1})\cdots
\bar s^{\ts t}_{p^{}_{n}p^{\ts\prime}_n}(q^{2n-2}u^{-1})&\\
{}\times s^{}_{p^{}_{n+1}p^{\ts\prime}_{n+1}}(q^{-2n}u)\cdots
s^{}_{p^{}_Np^{\ts\prime}_N}(q^{-2N+2}u)&,
\eal
\end{equation}
where the $\bar s^{\ts t}_{ij}(u)$ denote
the matrix elements of the transposed matrix $\overline S^{\tss t}(u)$.
\eth

\bex
If $N=2$ then
\beql{sdn2}
c(u)=\bar  s^{\ts t}_{11}(u^{-1})\ts s^{}_{22}(q^{-2}u)-\qin\ts
\bar s^{\ts t}_{21}(u^{-1})\ts s^{}_{12}(q^{-2}u).
\end{equation}
If $N=3$ then
\beql{sdn3}
\bal
c(u)&=\bar  s^{\ts t}_{22}(u^{-1})\ts s^{}_{11}(q^{-2}u)\ts s^{}_{33}(q^{-4}u)
+\bar  s^{\ts t}_{12}(u^{-1})\ts s^{}_{31}(q^{-2}u)\ts s^{}_{23}(q^{-4}u)\\
&{}+q^{-2}\ts\bar  s^{\ts t}_{21}(u^{-1})\ts s^{}_{32}(q^{-2}u)\ts s^{}_{13}(q^{-4}u)
-q\ts\bar  s^{\ts t}_{12}(u^{-1})\ts s^{}_{21}(q^{-2}u)\ts s^{}_{33}(q^{-4}u)\\
&-q^{-1}\ts\bar  s^{\ts t}_{32}(u^{-1})\ts s^{}_{11}(q^{-2}u)\ts s^{}_{23}(q^{-4}u)
-q^{-3}\ts\bar  s^{\ts t}_{31}(u^{-1})\ts s^{}_{22}(q^{-2}u)\ts s^{}_{13}(q^{-4}u).
\eal
\non
\end{equation}
\eex

\bpf[Proof of Theorem~\ref{thm:exp}]
Our proof is based on the properties of the quantum minors
of the matrices $T(u)$ and $\overline T(u)$.
Let $a_1,\dots,a_r$ and
$b_1,\dots,b_r$ be indices from the set $\{1,\dots,N\}$.
Introduce the elements
\beql{smin}
{s\ts}^{a_1\cdots\ts a_r}_{b_1\cdots\ts b_r}(u)=
\sum_{c_1<\cdots<c_r}
{t\ts}^{a_1\cdots\ts a_r}_{c_1\cdots\ts c_r}(u)\ts
{\bar t\ts}^{b_r\cdots\ts b_1}_{c_r\cdots\ts c_1}(q^{2r-2}\ts u^{-1}),
\end{equation}
where the indices $c_1,\dots,c_r$ run over the set $\{1,\dots,N\}$.
In particular, ${s\ts}^{a}_{b}(u)=s_{ab}(u)$, and
\beql{sdeteq}
{s\ts}^{1\cdots\ts N}_{1\cdots\ts N}(u)=c(u)
\end{equation}
by \eqref{sdetqdet}. We shall now derive a recurrent formula
for the elements \eqref{smin} with the conditions $b_i=a_i$
for $i=1,\dots,r-1$ and $a_1<\cdots<a_r$.
First, by the formulas \eqref{qminorgen}--\eqref{center}
for the quantum
minors we can write
\beql{sminp}
\bal
r!\ts{s\ts}^{a_1\cdots\ts a_r}_{a_1\cdots\ts a_{r-1},b_r}(u)&=
\sum_{c_1,\dots,c_r}
{t\ts}^{a_1\cdots\ts a_r}_{c_1\cdots\ts c_r}(u)\ts
{\bar t\ts}^{b_r,a_{r-1}\cdots\ts a_1}_{c_r\cdots\ts c_1}(q^{2r-2}\ts u^{-1})\\
{}&=r\ts\sum_{c_1,\dots,c_r}
{t\ts}^{a_1\cdots\ts a_r}_{c_1\cdots\ts c_r}(u)\ts
{\bar t\ts}^{a_{r-1}\cdots\ts a_1}_{c_{r-1}\cdots\ts c_1}(q^{2r-4}\ts u^{-1})
\ts \bar t^{}_{b_rc_r}(q^{2r-2}\ts u^{-1})\\
{}&=r\ts\sum_{c_1,\dots,c_r}
{\bar t\ts}^{a_{r-1}\cdots\ts a_1}_{c_{r-1}\cdots\ts c_1}(q^{2r-4}\ts u^{-1})
\ts{t\ts}^{a_1\cdots\ts a_r}_{c_1\cdots\ts c_r}(u)
\ts \bar t^{}_{b_rc_r}(q^{2r-2}\ts u^{-1}).
\eal
\end{equation}
Next, applying the column expansion to the quantum minor
${t\ts}^{a_1\cdots\ts a_r}_{c_1\cdots\ts c_r}(u)$, we obtain
\beql{sminpn}
\bal
{}&(r-1)!\ts{s\ts}^{a_1\cdots\ts a_r}_{a_1\cdots\ts a_{r-1},b_r}(u)\\
&=\sum_{c_1,\dots,c_r}\sum_{k=1}^r(-q)^{k-r}\ts
{\bar t\ts}^{a_{r-1}\cdots\ts a_1}_{c_{r-1}\cdots\ts c_1}(q^{2r-4}\ts u^{-1})
\ts{t\ts}^{a_1\cdots\wh{a}_k\cdots\ts a_r}_{c_1\cdots\ts c_{r-1}}(u)
\ts t^{}_{a_kc_r}(q^{-2r+2}\ts u)\ts\bar t^{}_{b_rc_r}(q^{2r-2}\ts u^{-1})\\
&=\sum_{c_1,\dots,c_{r-1}}\sum_{k=1}^r(-q)^{k-r}\ts
{\bar t\ts}^{a_{r-1}\cdots\ts a_1}_{c_{r-1}\cdots\ts c_1}(q^{2r-4}\ts u^{-1})
\ts{t\ts}^{a_1\cdots\wh{a}_k\cdots\ts a_r}_{c_1\cdots\ts c_{r-1}}(u)
\ts s^{}_{a_kb_r}(q^{-2r+2}\ts u),
\eal
\non
\end{equation}
where the hats indicate the symbols to be omitted.
Now for any $k<r$ write
\beql{tk}
{\bar t\ts}^{a_{r-1}\cdots\ts a_1}_{c_{r-1}\cdots\ts c_1}(q^{2r-4}\ts u^{-1})
=(-q)^{k-1}\ts
{\bar t\ts}^{a_{r-1}\cdots\wh{a}_k\cdots\ts a_1,a_k}_{c_{r-1}\cdots\ts c_1}
(q^{2r-4}\ts u^{-1})
\end{equation}
and apply the column expansion to this minor to
bring the previous formula to the form
\beql{sminpnf}
\bal
{}&(r-2)!\ts{s\ts}^{a_1\cdots\ts a_r}_{a_1\cdots\ts a_{r-1},b_r}(u)\\
&=\sum_{c_1,\dots,c_{r-1}}\Big\{
(-q)^{r-2}\ts\bar t^{}_{a_{r-1}c_1}(u^{-1})
\ts{\bar t\ts}^{a_{r-2}\cdots\ts a_1}_{c_{r-1}\cdots\ts c_2}
(q^{2r-4}\ts u^{-1})
\ts{t\ts}^{a_1\cdots\ts a_{r-1}}_{c_1\cdots\ts c_{r-1}}(u)
\ts s^{}_{a_rb_r}(q^{-2r+2}\ts u)\\
{}&+\sum_{k=1}^{r-1}(-q)^{2k-r-1}\ts\bar t^{}_{a_kc_1}(u^{-1})
\ts{\bar t\ts}^{a_{r-1}\cdots\wh{a}_k\cdots\ts a_1}_{c_{r-1}\cdots\ts c_2}
(q^{2r-4}\ts u^{-1})
\ts{t\ts}^{a_1\cdots\wh{a}_k\cdots\ts a_r}_{c_1\cdots\ts c_{r-1}}(u)
\ts s^{}_{a_kb_r}(q^{-2r+2}\ts u)\Big\}.
\eal
\non
\end{equation}
Applying again \eqref{center}, we write this as
\beql{sminpnfc}
\bal
{}&(r-2)!\ts{s\ts}^{a_1\cdots\ts a_r}_{a_1\cdots\ts a_{r-1},b_r}(u)\\
&=\sum_{c_1,\dots,c_{r-1}}\Big\{
(-q)^{r-2}\ts\bar t^{}_{a_{r-1}c_1}(u^{-1})
\ts{t\ts}^{a_1\cdots\ts a_{r-1}}_{c_1\cdots\ts c_{r-1}}(u)
\ts{\bar t\ts}^{a_{r-2}\cdots\ts a_1}_{c_{r-1}\cdots\ts c_2}
(q^{2r-4}\ts u^{-1})\ts s^{}_{a_rb_r}(q^{-2r+2}\ts u)\\
{}&+\sum_{k=1}^{r-1}(-q)^{2k-r-1}\ts\bar t^{}_{a_kc_1}(u^{-1})
\ts{t\ts}^{a_1\cdots\wh{a}_k\cdots\ts a_r}_{c_1\cdots\ts c_{r-1}}(u)
\ts{\bar t\ts}^{a_{r-1}\cdots\wh{a}_k\cdots\ts a_1}_{c_{r-1}\cdots\ts c_2}
(q^{2r-4}\ts u^{-1})\ts s^{}_{a_kb_r}(q^{-2r+2}\ts u)\Big\}.
\eal
\non
\end{equation}
Finally, using the column expansion \eqref{columnexp}
and the definition \eqref{barSuorth} of the elements $\bar s_{ij}(u)$
we get the following recurrence relation
\beql{recs}
\bal
{s\ts}^{a_1\cdots\ts a_r}_{a_1\cdots\ts a_{r-1},b_r}(u)
&={\bar s\ts}^{t}_{a_{r-1}a_{r-1}}(u^{-1})\ts
{s\ts}^{a_1\cdots\ts a_{r-2}}_{a_1\cdots\ts a_{r-2}}(q^{-2}\ts u)
\ts s^{}_{a_rb_r}(q^{-2r+2}\ts u)\\
&+\sum_{l=1}^{r-2}(-q)^{2r-2l-3}\ts{\bar s\ts}^{t}_{a_{l}a_{r-1}}(u^{-1})\ts
{s\ts}^{a_1\cdots\wh{a}_l\cdots\ts a_{r-1}}_{a_1
\cdots\wh{a}_l\cdots\ts a_{r-2},a_l}(q^{-2}\ts u)
\ts s^{}_{a_rb_r}(q^{-2r+2}\ts u)\\
&+\sum_{k=1}^{r-1}\Big\{(-q)^{2k-2r+1}\ts{\bar s\ts}^{t}_{a_{r}a_{k}}(u^{-1})
\ts {s\ts}^{a_1\cdots\wh{a}_k\cdots\ts a_{r-1}}_{a_1
\cdots\wh{a}_k\cdots\ts a_{r-1}}(q^{-2}\ts u)
\ts s^{}_{a_kb_r}(q^{-2r+2}\ts u)\\
&\quad+\sum_{l=1}^{k-1}(-q)^{2k-2l-2}\ts{\bar s\ts}^{t}_{a_{l}a_{k}}(u^{-1})
\ts {s\ts}^{a_1\cdots\wh{a}_l\cdots\wh{a}_k\cdots\ts a_{r}}_{a_1
\cdots\wh{a}_l\cdots\wh{a}_k\cdots\ts a_{r-1},a_l}(q^{-2}\ts u)
\ts s^{}_{a_kb_r}(q^{-2r+2}\ts u)\\
&\quad+\sum_{l=k+1}^{r-1}(-q)^{2k-2l}\ts{\bar s\ts}^{t}_{a_{l}a_{k}}(u^{-1})
\ts {s\ts}^{a_1\cdots\wh{a}_k\cdots\wh{a}_l\cdots\ts a_{r}}_{a_1
\cdots\wh{a}_k\cdots\wh{a}_l\cdots\ts a_{r-1},a_l}(q^{-2}\ts u)
\ts s^{}_{a_kb_r}(q^{-2r+2}\ts u)\Big\}.
\eal
\non
\end{equation}
Using \eqref{sdeteq}
and starting with ${s\ts}^{1\cdots\ts N}_{1\cdots\ts N}(u)$
we apply this recurrence relation repeatedly to get
an explicit expression for the series $c(u)$ in terms
of the generators $s_{ij}(u)$ and $\bar s_{ij}(u)$.
It follows from the definition of the map $\pi^{}_N$ that
$c(u)$ will be written as a combination
of the monomials of the required form, the coefficients
being powers of $-q$. The exact values of the powers
are easily found by calculating the number of inversions
of the permutations occurring in the recurrence relation.
\epf

\bre\label{rem:twy}
The argument used in the proof of Theorem~\ref{thm:exp}
can also be applied to produce a simpler proof for
the formula for the Sklyanin determinant
given in \cite{m:ya} in the
case of orthogonal
and symplectic twisted Yangians. However, the argument does not seem to be
directly applicable
to the case of the symplectic twisted $q$-Yangian $\Y^{\rm tw}_q(\spa_{2n})$.
\ere

The image of the matrix $\overline S(u)$ under
the evaluation homomorphism \eqref{evaluaff} is found from
the relation \eqref{relsbarsorth} so that
\beql{sbarim}
\overline S(u)\mapsto \frac{1+uq^{-1}}{1+uq}(\overline S+q\ts u \ts S).
\end{equation}

Applying the evaluation homomorphism
to the Sklyanin determinant $\sdet S(u)$
and using \eqref{evaluaff} and \eqref{sbarim} we derive
the following corollary from Theorem~\ref{thm:exp}.

\bco\label{cor: sdetfin}
The coefficients of the polynomial
\beql{sdetfin}
\bal
C(u)=\sum_{p\in \Sym_{N}} (-q)^{-l(p)+l(p')}
\big[u\ts \overline S+q\ts S\big]_{p^{\ts\prime}_1p^{}_1}\cdots
\big[u\ts \overline S+q^{2n-1}S\big]&_{p^{\ts\prime}_np^{}_{n}}\\
{}\times \big[u\ts S+q^{2n-1}\ts
\overline S\ts\big]_{p^{}_{n+1}p^{\ts\prime}_{n+1}}\cdots
\big[u\ts S+q^{2N-3}\ts\overline S\ts\big]&_{p^{}_Np^{\ts\prime}_N}
\eal
\non
\end{equation}
are Casimir elements for the algebra $\U^{\rm tw}_q(\oa_{N})$.
Moreover, the polynomial $C(u)$ is monic of degree $N$.
\eco

\bpf
The polynomial $C(u)$ is obtained by the application of
the evaluation homomorphism
to the series $c(u)$
and multiplication by an appropriate rational function in $u$.
The centrality of its coefficients thus follows from
Theorem~\ref{thm:exp}. Obviously, the degree of $C(u)$
does not exceed $N$. The coefficient of $u^N$
can only occur in the summands with the property $p=p'$.
However, it follows from the definition
of the map \eqref{mapp} that this property
is satisfied by the only permutation $p$
\beql{onlyp}
p=\begin{cases}
        (N-1,N-3,\dots,1,2,4,\dots,N)\qquad
        &\text{if}\quad N\quad\text{is even,}\\
        (N-1,N-3,\dots,2,1,3,\dots,N)\qquad
        &\text{if}\quad N\quad\text{is odd.}
   \end{cases}
\end{equation}
Since the diagonal entries of $S$ and $\overline S$ are ones,
the second statement follows.
\epf

The polynomial $C(u)$ may be regarded as a $q$-analog of the Capelli
polynomial for the algebra $\U^{\rm tw}_q(\oa_{N})$. It would be interesting
to find its eigenvalues in the irreducible representations
and to get the corresponding $q$-analogs of the Capelli identities;
cf. \cite{m:ya}, \cite{gi:ce}, \cite{hkp:ce}, \cite{nuw:dp}.

\bex\label{ex:o3}
If $N=3$ then
\beql{o3}
C(u)=(u+q)\big[(u+q)(u+q^3)-q^{2}u\ts C\big],
\end{equation}
where
\beql{casim}
C=s^{2}_{21}+q^2\ts s^{2}_{32}+s^{2}_{31}-q\ts s^{}_{21}\ts
s^{}_{32}\ts s^{}_{31}.
\end{equation}
\eex

The following corollary
provides a
characteristic identity for the algebra $\U^{\rm tw}_q(\oa_{N})$;
cf. \cite{m:ya}, \cite{nt:yg}.

\bco\label{cor:charid}
We have the identity
\beql{charid}
C(-q^{2N-3}\ts\overline S\ts S^{-1})=0.
\end{equation}
\eco

\bpf
Introduce the quantum comatrix $\wh S(u)$ by the formula
\beql{comatr}
\wh S(u)\ts S(uq^{-2N+2})=\sdet S(u).
\end{equation}
Explicit expressions for the matrix elements of the matrix
$\wh S(u)$ can be found from the recurrence relation
for the elements ${s\ts}^{a_1\cdots\ts a_r}_{a_1\cdots\ts a_{r-1},b_r}(u)$;
see the end of the proof of Theorem~\ref{thm:exp}.
Using these expressions and
applying the evaluation homomorphism to both sides of
\eqref{comatr} we get
\beql{comatrc}
C(u)=\wh C(u)\ts (\overline S+uq^{-2N+3}\ts S),
\end{equation}
where $\wh C(u)$ is a polynomial in $u$ with coefficients
in the algebra $\U^{\rm tw}_q(\oa_{N})\ot\End \C^N$. This completes the proof.
\epf

\end{document}